\documentclass[12pt]{amsart}
\usepackage{amsmath, amssymb, amsthm, bm, float, mathtools, enumitem, caption}
\usepackage{hyperref}

\usepackage[noabbrev,capitalize]{cleveref}
\usepackage{adjustbox}
\crefname{equation}{}{}
\usepackage{ytableau}
\usepackage{graphicx, tikz}
\usepackage[textheight=8.5in, textwidth=6.75in]{geometry}
\usepackage{subcaption}
\usepackage{microtype}
\usepackage{stmaryrd}
\usepackage{mathtools}

\newtheorem{theorem}{Theorem}[section]
\newtheorem{lemma}[theorem]{Lemma}
\newtheorem{corollary}[theorem]{Corollary}
\newtheorem{proposition}[theorem]{Proposition}
\newtheorem{conjecture}[theorem]{Conjecture}
\newtheorem*{conjecture*}{Conjecture}

\theoremstyle{definition}

\theoremstyle{remark}

\numberwithin{equation}{section}

\newcommand*\diff{\mathop{}\!\mathrm{d}}

\title[The parity of the number of missing integers in partitions]{The parity of the number of missing integers in partitions}
\date{\today}
\author{Julia Q.D. Du \and Gurinder Singh}
\address{School of Mathematical Sciences, Hebei Key Laboratory of Computational Mathematics and Applications, Hebei Research Center of the Basic
discipline Pure Mathematics, Hebei Normal University, Shijiazhuang 050024, P.R. China}
\email{qddu@hebtu.edu.cn}
\address{Postdoctoral Research Station of Mathematics, Hebei Normal University, Shijiazhuang 050024, P.R. China}
\email{gurindermaan1018@gmail.com}

\begin{document}

\begin{abstract}
Missing integers in partitions have been studied by many authors.
A missing integer of a partition is a positive integer smaller than its largest part that does not occur as a part. We investigate parity biases for the number of missing integers in partitions and overpartitions. For every integer $c\ge2$, we prove that the coefficients $a_c(n)$ of $(-cq;q)_\infty/(-q;q)_\infty$ are eventually strictly positive and provide an effectively computable threshold. The proof uses the circle method, combining a saddle-point analysis on the major arc with explicit uniform estimates on the minor arcs. Specializing to $c=2$ and $c=3$, and verifying the remaining finite cases, we show that partitions and overpartitions of $n$ with an even number of missing integers outnumber those with an odd number for $n>34$ and $n>27$, respectively. We further prove that, for $n>10$, the difference between the numbers of ordinary partitions with an even and an odd number of missing integers is smaller than the number of partitions of $n$ into distinct parts. These results confirm three conjectures of Bhoria, Eyyunni and Santra. We also give a combinatorial interpretation of $a_c(n)$ for every integer $c\ge4$ in terms of colored partitions, generalizing the partition and overpartition interpretations of Bhoria, Eyyunni and Santra for $c=2$ and $c=3$, respectively.
\end{abstract}

\maketitle
\section{Introduction}
Partition statistics provide a natural framework for investigating the interplay among combinatorial structure, generating functions, and arithmetic properties.
Classical statistics, including the number of
parts, the largest part, the rank, and the crank,
depend on the parts occurring in a partition.
A complementary perspective is
obtained by studying the positive integers that are absent.
Grabner and Knopfmacher \cite{GK_2006}
investigated gaps in partitions and
showed that almost all partitions contain gaps, thereby motivating the
study of statistics associated with missing integers.
\par One such statistic is the \emph{minimal excludant}, or mex, of a
partition: the smallest positive integer that is not a part of the partition.
Grabner and Knopfmacher \cite{GK_2006} studied this statistic under the
name ``the least gap''.
Andrews and Newman
\cite{AN_2019} subsequently
developed its combinatorial theory using the terminology
``minimal excludant''.
They established the following elegant identity involving the quantity 
$\sigma$mex($n$), which denotes the sum of the
minimal excludants over all partitions of $n$,
\begin{align}\label{eq-AN}
\sum_{n=0}^{\infty} \sigma\mathrm{mex}(n)q^n = (-q;q)_{\infty}^2 = \sum_{n=0}^{\infty} \mathcal{D}_2(n)q^n,
\end{align}
where $\mathcal{D}_2(n)$ counts the number of partitions of 
$n$ into distinct parts with two colors.
Here and throughout, we use the standard notation
\[(a;q)_{\infty}:=\prod_{j=0}^{\infty}(1-aq^j), \qquad |q|<1.\]
\par The identity \eqref{eq-AN}
had already appeared in the work of Grabner and Knopfmacher \cite{GK_2006};
a combinatorial proof was later given by Ballantine and Merca \cite{BM_2021}. 
Further connections between mex and
classical partition statistics were established by Hopkins, Sellers,
and Stanton \cite{HSS_2022}, who investigated its relationship with the
crank. Baruah et al. \cite{BEM_2023} refined the arithmetic function $\sigma\mathrm{mex}(n)$ by considering the sums of odd and even minimal excludants separately. 
For more information, please refer to \cite{BS_2021, CR_2021, SS_2020}.
\par Another related statistic is the \emph{maximal excludant}, or maex, introduced
by Chern \cite{Ch_2021},
that is, the largest nonnegative integer less than the largest part
that does not occur as a part of the partition.
Whenever a positive integer below the largest part is
missing, mex and maex identify the smallest and largest such integers,
respectively. These statistics therefore capture the extremal elements
of the set of missing integers.
\par Bhoria, Eyyunni and Santra \cite{BES_2026} 
shifted attention from these
extremal values to the entire set of missing positive integers less than
the largest part  of a partition. For example, the partition $6+4+4+1$ has three missing integers, $5$, $3$, and $2$.
In \cite{BES_2026}, the authors studied the number of missing integers in
partitions and overpartitions, obtaining generating functions,
identities, and congruences.
Aniceto and Ballantine \cite{AB_2026}  subsequently provided combinatorial proofs
for several of their results.
\par Let $\mathcal{M}_e(n)$ (resp., $\mathcal{M}_o(n)$) denote the number of partitions of a positive integer $n$ with an even (resp., odd) number of missing integers less than the largest part. Set $\mathcal{M}_e(0)=1$ and $\mathcal{M}_o(0)=0$. For example, the partitions $5$, $4+1$, $2+2+1$, $2+1+1+1$, and $1+1+1+1+1$ of $5$ have an even number of missing integers, whereas $3+2$ and $3+1+1$ have an odd number of missing integers. Therefore, $\mathcal{M}_e(5)=5$ and $\mathcal{M}_o(5)=2$.
Bhoria, Eyyunni and Santra \cite{BES_2026} proved that
\begin{align*}
\sum_{n=0}^{\infty}\left(\mathcal{M}_e(n)-\mathcal{M}_o(n)\right)q^n=\frac{(-2q;q)_{\infty}}{(-q;q)_{\infty}}.
\end{align*}
The coefficients of this product are recorded as the OEIS sequence A268498 \cite{S_onlie}.
Based on the apparent eventual positivity
of these coefficients, they proposed a conjecture asserting that
partitions with an even number of missing integers outnumber those
with an odd number whenever $n>34$ (see \cite[Conjecture 2.7]{BES_2026}).
\begin{conjecture}{\cite[Conjecture 2.7]{BES_2026}}\label{conj1}
For $n>34$, we have
\begin{align*}
\mathcal{M}_e(n) >\mathcal{M}_o(n).
\end{align*}
\end{conjecture}
An overpartition of a positive integer $n$ is a partition of $n$ in which the first occurrence of a part may be overlined. 
Let $\overline{\mathcal{M}}_e(n)$ (resp., $\overline{\mathcal{M}}_o(n)$) denote the number of overpartitions of a positive integer $n$ with an even (resp., odd) number of missing integers less than the largest part. Set $\overline{\mathcal{M}}_e(0)=1$ and $\overline{\mathcal{M}}_o(0)=0$.
Bhoria, Eyyunni and Santra \cite{BES_2026} proved that
\begin{align*}
\sum_{n=0}^{\infty}\left(\overline{\mathcal{M}}_e(n)-\overline{\mathcal{M}}_o(n)\right)q^n=\frac{(-3q;q)_{\infty}}{(-q;q)_{\infty}}.
\end{align*}
The authors of \cite{BES_2026} also proposed the following overpartition analogue of Conjecture \ref{conj1}.
\begin{conjecture}{\cite[Conjecture 3.6]{BES_2026}}\label{conj2}
For $n>27$, we have
\begin{align*}
\overline{\mathcal{M}}_e(n) >\overline{\mathcal{M}}_o(n).
\end{align*}
\end{conjecture}
In this article, we prove Conjectures \ref{conj1} and \ref{conj2}, and establish a common generalization for the family 
\begin{align}
F_c(q):=\sum_{n=0}^{\infty}a_c(n)q^n=\frac{(-cq;q)_{\infty}}{(-q;q)_{\infty}},
\label{gf-Fc}
\end{align}
where $c\ge2$ is an integer. The obvious question that arises here is the following. Can we give a combinatorial interpretation to $a_c(n)$ for $c\geq4$? The answer generalizes the partition definitions of $a_c(n)$ for $c=2,3$. 
\par Let $c\geq2$. We construct colored partitions by coloring each distinct part size of a partition independently with one of $c-1$ colors. All occurrences of the same part size receive the same color and different sizes may receive the same color. Let $\mathcal{M}^c_e(n)$ and $\mathcal{M}^c_o(n)$ be the numbers of these colored partitions of a positive integer $n$ having, respectively, an even and an odd number of missing integers. Set $\mathcal{M}^c_e(0)=1$ and $\mathcal{M}^c_o(0)=0$.
\begin{theorem}\label{Thm_CI}
Let $c\geq2$ be an integer. For all $n\geq0$,
\begin{align*}
a_c(n)=\mathcal{M}^c_e(n)-\mathcal{M}^c_o(n).
\end{align*}
\end{theorem}
For $c=2$, there is only one color, so these become the corresponding ordinary partition counts, i.e., $\mathcal{M}^2_e(n)=\mathcal{M}_e(n)$ and $\mathcal{M}^2_o(n)=\mathcal{M}_o(n)$. For $c=3$, there are two available colors, say red ($R$) and blue ($B$). For example, partition $(2,1,1)$ has four possible colorings: $(2_R,1_R,1_R)$, $(2_R,1_B,1_B)$, $(2_B,1_R,1_R)$, and $(2_B,1_B,1_B)$. Replace every red part with a nonoverlined part of the same size. For each blue part size, overline its first occurrence and leave its remaining occurrences nonoverlined. Then, we obtain $\mathcal{M}^3_e(n)=\overline{\mathcal{M}}_e(n)$ and $\mathcal{M}^3_o(n)=\overline{\mathcal{M}}_o(n)$. Thus, Theorem \ref{Thm_CI} is a generalization of the combinatorial interpretations of $a_c(n)$ for $c=2,3$.
\par For every $c\geq2$, we show that the coefficients of $F_c(q)$ are strictly positive for all sufficiently large indices and provide an effectively computable threshold.
\begin{theorem}\label{Thm_main}
For every integer $c\geq2$, there exists a positive integer $N_c$ such that for all $n\geq N_c$,
\begin{align*}
a_c(n)>0.
\end{align*}
\end{theorem}
Specializing Theorem \ref{Thm_main} to $c=2$ and $c=3$, and verifying the remaining finite cases, we prove Conjectures \ref{conj1} and \ref{conj2}.
\begin{corollary}\label{cor1}
Conjectures \ref{conj1} and \ref{conj2} are true. 
\end{corollary}
We prove Theorem \ref{Thm_main} using the circle method with explicit error bounds, a strategy also employed by Chern \cite{Ch_2023} and Craig \cite{Craig_2023} in establishing positivity results for coefficients of nonmodular $q$-products. For every integer $c\ge2$, the generating function $F_c(q)$ is nonmodular, which complicates the proof of eventual positivity of its coefficients. In particular, the modular transformation laws underlying the usual Rademacher-type coefficient formulas cannot be applied directly to $F_c(q)$. We overcome this difficulty by deriving explicit estimates that isolate a positive major arc contribution and show that it dominates the contributions from the minor arcs for $n\ge N_c$, thereby establishing eventual positivity with an effective threshold. On the major arc near $q=1$, Euler–Maclaurin summation and a saddle-point analysis yield a positive lower bound for the coefficient integral. On the minor arcs, Jacobi’s triple-product identity, Poisson summation, and the transformation law for Dedekind’s eta-function provide uniform upper bounds near roots of unity.
\par Let $p(n)$ denote the number of ordinary partitions of $n$, and $q(n)$ denote the number of partitions of $n$ into distinct parts. In \cite{BES_2026}, Bhoria, Eyyunni and Santra conjectured the following upper bound for the difference $\mathcal{M}_e(n) -\mathcal{M}_o(n)$.
\begin{conjecture}{\cite[Conjecture 2.8]{BES_2026}}\label{conj3}
For $n>10$, we have
\begin{align*}
q(n)>\mathcal{M}_e(n) -\mathcal{M}_o(n).
\end{align*}
\end{conjecture}
In this article, we also prove Conjecture \ref{conj3}.
\begin{theorem}\label{Thm_main_2}
Conjecture \ref{conj3} is true.
\end{theorem}
The proof of Theorem \ref{Thm_main_2} combines explicit analytic bounds for $a_2(n)=\mathcal{M}_e(n) -\mathcal{M}_o(n)$ and $q(n)$. We obtain an exponential upper bound for $a_2(n)$ by estimating $F_2(e^{-t})$ and using the eventual positivity established in the proof of Theorem \ref{Thm_main}. For the lower bound, an explicit injection gives $q(n)\ge p(\lfloor n/2\rfloor)$, to which we apply the estimate of Banerjee, Paule, Radu and Zeng \cite[Theorem 1.3]{BPRZ} for $p(n)$. Then we compare the resulting bounds to prove $q(n)>a_2(n)$ for $n\ge2500$, and the remaining cases $11\le n\le2499$ are verified computationally.
\par The rest of this paper is organized as follows. Section \ref{Sec2} collects the necessary facts about dilogarithms, Euler–Maclaurin summation, the Dedekind eta-function, and Jacobi’s triple-product identity. We prove Theorem \ref{Thm_CI} in Section \ref{Sec3new}. In Section \ref{Sec3}, we prove Theorem \ref{Thm_main} by comparing the major and minor arc contributions, obtain an effective threshold for coefficient positivity, and establish Corollary \ref{cor1} by verifying the remaining finite cases. Section \ref{Sec4} proves Theorem \ref{Thm_main_2} through explicit upper and lower bounds for the parity difference and the distinct partition function, respectively. Section \ref{Sec5} presents concluding remarks.
\section{Preliminaries}\label{Sec2}
In this section, we collect the analytic tools used in the proof of Theorem \ref{Thm_main}. These include basic properties of the dilogarithm function, Euler-Maclaurin summation with an explicit remainder, the transformation law for Dedekind's $\eta$-function, and Jacobi’s triple-product identity.
\par First, we fix some notation. For a complex number $z$, let $\mathrm{Re}(z)$ and $\mathrm{Im}(z)$ denote its real and imaginary parts, respectively. Let $\mathbb{H}$ denote the upper half-plane, that is,
$\mathbb{H}:=\{z\in\mathbb{C}:\mathrm{Im}(z)>0\}$. 
We use $\operatorname{Log}(\cdot)$ for the principal complex logarithm. The notation $\log(\cdot)$ denotes the real natural logarithm or an explicitly specified analytic logarithm.
\vspace{0.1cm}

\textbf{The dilogarithm.} 
The dilogarithm function is defined by
\begin{align*}
\text{Li}_2(z)=\sum_{m=1}^{\infty}\frac{z^m}{m^2}, \qquad|z|<1.
\end{align*}
Its principal branch extends analytically to $\mathbb C\setminus[1,\infty)$ and has the integral representation
\begin{align*}
\text{Li}_2(z)=-\int_{0}^{z}\frac{\operatorname{Log}(1-u)}{u}\diff u.
\end{align*}
The dilogarithm function has the following properties:
\begin{enumerate}
    \item For any $z\in\mathbb C\setminus[0,\infty)$, we have
    \begin{align}\label{rel-Li2}
\text{Li}_2\left(\frac{1}{z}\right)+\text{Li}_2(z)=-\frac{\pi^2}{6}-\frac{1}{2}\operatorname{Log}^2(-z);
\end{align}

\item 
The following duplication identity holds:
\begin{align}\label{Li_id}
\text{Li}_2(z)+\text{Li}_2(-z)=\frac{\text{Li}_2(z^2)}{2},\qquad |z|<1;
\end{align}

\item 
For every real $z>0$, we have $-z<\text{Li}_2(-z)<0$;

\item 
Values for some specific $z$ are known, for example,
$\text{Li}_2(-1)=-\pi^2/12$, and
$\text{Li}_2(1/2)=\pi^2/12-(\log2)^2/2$.

\end{enumerate}
For more on this function, please see \cite{Zagier}.
\vspace{0.1cm}

\textbf{The Dedekind $\eta$-function.} 
We first recall the definition of the full modular group, that is,
\begin{align*}
\text{SL}_2(\mathbb{Z})=\bigg\{\!\begin{pmatrix}a &b\\ c &d\end{pmatrix}\colon a,b,c,d\in\mathbb{Z},
\;\;\textup{and}\;\;ad-bc=1\bigg\}.
\end{align*}
The matrix $\gamma=\begin{pmatrix}a &b\\ c &d\end{pmatrix}$ acts on
$\tau\in\mathbb{H}$ by the linear fractional transformation
\begin{align*}
\gamma\tau=\dfrac{a\tau+b}{c\tau+d}.
\end{align*}

The Dedekind eta function is given by
\begin{align*}
\eta(\tau)=e^{\pi i\tau/12}\prod_{j=1}^{\infty}\left(1-e^{2\pi ij\tau}\right),\ \ \tau\in\mathbb{H}.
\end{align*}
For any $\tau\in\mathbb{H}$ and $\gamma=\begin{pmatrix}
a & b \\ c & d 
\end{pmatrix}\in \text{SL}_2(\mathbb{Z})$, with $c>0$,
the transformation formula \cite[Theorem 3.4]{Apostal1990} for $\eta(\tau)$ is given by 
\begin{align}\label{eq_eta_TF}
\eta(\gamma\tau)=\zeta_\gamma(c\tau+d)^{1/2}\eta(\tau),
\end{align}
where $|\zeta_\gamma|=1$.

\vspace{0.1cm}

\textbf{Euler–Maclaurin summation.} We use Euler–Maclaurin summation with an explicit integral remainder; see \cite[pp. 279-292]{Olver}. Let \(g:[0,\infty)\to\mathbb C\) be twice continuously differentiable. Suppose that \(g\) and \(g''\) are absolutely integrable,  \(g(u),g'(u)\to0\) as \(u\to\infty\), and \(\sum_{m=1}^\infty g(m)\) converges. Then
\begin{align}
\sum_{m=1}^\infty g(m)
=
\int_0^\infty g(u)\,du
-\frac{g(0)}2
-\frac{g'(0)}{12}
+\mathcal{R}_2(g), \label{eq-g-EM-sum}
\end{align}
where
\begin{align*}
\mathcal{R}_2(g)=-\frac{1}{2}\int_0^\infty B_2(u)\,g''(u)\diff u,
\end{align*}
and
\begin{align*}
B_2(u)=\{u\}^2-\{u\}+\frac{1}{6}
\end{align*}
is the second Bernoulli function and $\{u\}$ is the fractional part of $u$.
Notice that $B_2(u)$ satisfies $\sup_{u\geq 0}|B_2(u)|=\frac{1}{6}$, so
\begin{align*}
\vert\mathcal{R}_2(g)\vert\le\frac1{12}\int_0^\infty|g''(u)|\diff u.
\end{align*}
The complex-valued version follows by applying the real-valued formula separately to the real and imaginary parts. 


\textbf{Jacobi's triple-product identity.}
For $|q|<1$ and $z\ne0$, Jacobi's triple-product identity \cite[Theorem 1.3.3]{Berndt} states that
\begin{align}\label{JTPI*}
\sum_{m=-\infty}^\infty z^m q^{m^2}
=(-zq;q^2)_\infty(-q/z;q^2)_\infty(q^2;q^2)_\infty.
\end{align}
Replacing $z$ by $-z/q$ in \eqref{JTPI*}, and then $q$ by $q^{1/2}$, we obtain
\begin{align}\label{JTPI}
\sum_{m=-\infty}^\infty (-z)^m q^{m(m-1)/2}
=(q;q)_\infty(z;q)_\infty(q/z;q)_\infty.
\end{align}
\section{Proof of Theorem \ref{Thm_CI}}\label{Sec3new}
Let $\mathcal{P}(n)$ denote the set of partitions of $n$. For a nonempty partition
$\lambda=(\lambda_1,\ldots,\lambda_{r})$, with
$\lambda_1\geq\cdots\geq\lambda_{r}>0$, we write $\ell(\lambda)$ for its number of parts, counted with multiplicity; $d(\lambda)$ for its number of distinct part sizes; $M(\lambda)$ for the number of positive integers smaller than $\lambda_1$
that do not occur as parts of $\lambda$. The conjugate of $\lambda$ is defined by $\lambda'=(\lambda_1',\lambda_2',\ldots,\lambda_s')$, where $\lambda_j'$ is the number of parts of $\lambda$ that are greater than or equal to $j$ (see \cite[p. 7]{And_1976}).
Since precisely $d(\lambda)-1$ distinct part sizes are smaller than $\lambda_1$,
we have
\begin{equation}\label{eq:m1}
M(\lambda)=(\lambda_1-1)-(d(\lambda)-1)
=\lambda_1-d(\lambda).
\end{equation}

\begin{proposition}\label{prop:main}
For $c\geq 2$ and $n\geq1$,
\begin{align}
a_c(n)=\sum_{\lambda\in\mathcal{P}(n)}
(-1)^{M(\lambda)}(c-1)^{d(\lambda)}.
\label{eq:m2}
\end{align}
Also, $a_c(0)=1$, corresponding to the empty partition.
\end{proposition}

\begin{proof}
For each possible part size $j\geq1$, expand its factor in the infinite product on the R.H.S. of \eqref{gf-Fc} as
\begin{equation}\label{eq:l1}
\frac{1+cq^j}{1+q^j}
=1+(c-1)\frac{q^j}{1+q^j}
=1+(c-1)\sum_{m=1}^{\infty}(-1)^{m-1}q^{jm}.
\end{equation}
The constant term corresponds to omitting the part size $j$. Choosing the term
with index $m\geq1$ corresponds to including exactly $m$ copies of $j$, with
weight $(c-1)(-1)^{m-1}$.

Write a partition in multiplicity notation as
$\lambda=(1^{m_1},2^{m_2},\ldots)$. The product of its local weights is
\begin{align*}
\prod_{j:m_j>0}(c-1)(-1)^{m_j-1}
&=(c-1)^{d(\lambda)}
(-1)^{\sum_{j:m_j>0}(m_j-1)}=(c-1)^{d(\lambda)}
(-1)^{\ell(\lambda)-d(\lambda)}.
\end{align*}
Its power of $q$ is $\sum_{j\geq1}jm_j=|\lambda|$. Multiplying
\eqref{eq:l1} over $j$ and taking the coefficient of $q^n$ gives
\begin{align}\label{eq:l2}
a_c(n)=\sum_{\lambda\in\mathcal{P}(n)}
(-1)^{\ell(\lambda)-d(\lambda)}(c-1)^{d(\lambda)}.
\end{align}
\par Let $\mu=\lambda'$.
Then $|\mu|=|\lambda|$ and $\mu_1=\ell(\lambda)$. Moreover, $d(\mu)=d(\lambda)$.
To see this, set $\lambda_{r+1}=0$. The multiplicity of the part size $i$ in $\mu$ is $\lambda_i-\lambda_{i+1}$. Thus the number of distinct
part sizes in $\mu$ equals the number of indices $i$ for which
$\lambda_i>\lambda_{i+1}$, which is exactly $d(\lambda)$.
Therefore, by \eqref{eq:m1},
\[
M(\mu)=\mu_1-d(\mu)
=\ell(\lambda)-d(\lambda).
\]
Conjugation is a bijection on $\mathcal{P}(n)$ satisfying
$d(\mu)=d(\lambda)$, so changing variables
from $\lambda$ to $\mu=\lambda'$ in \eqref{eq:l2} proves
\eqref{eq:m2}.
\end{proof}
Suppose now that $c\geq2$ is an integer. Color each distinct part size of a
partition independently with one of $c-1$ colors. All occurrences of the same
part size receive the same color and different sizes may receive the same color. Thus, a partition $\lambda$ admits exactly
$(c-1)^{d(\lambda)}$ colorings. Missing integers are determined by the underlying
partition, independently of its colors. Therefore,
\begin{align*}
\mathcal{M}^c_e(n)=\sum_{\substack{\lambda\in\mathcal{P}(n)\\ M(\lambda)\text{ even}}}
(c-1)^{d(\lambda)},\qquad
\mathcal{M}^c_o(n)=\sum_{\substack{\lambda\in\mathcal{P}(n)\\ M(\lambda)\text{ odd}}}
(c-1)^{d(\lambda)}.
\end{align*}
Thus, Proposition \ref{prop:main} proves Theorem \ref{Thm_CI}.

\section{Proof of Theorem \ref{Thm_main}}\label{Sec3}
Fix an integer $c\ge2$, and let $x>0$ be sufficiently small. By Cauchy's integral formula,
\begin{align*}
a_c(n)=\frac{1}{2\pi i}\int_{C}\frac{F_c(q)}{q^{n+1}} \diff q,
\end{align*}
where $C$ is the positively oriented circle $|q|=e^{-x}$. 
Let $q=e^{-z}$ with $z=x-iy$ and $-\pi\le y\le\pi$.
Then $\diff q=iq\diff y$, so
\begin{align}\label{eq3.1}
a_c(n)=\frac{e^{nx}}{2\pi }\int_{-\pi}^{\pi}F_c(e^{-x+iy})e^{-iny} \diff y.
\end{align}
The integral in \eqref{eq3.1} is real by conjugation symmetry. For brevity, let $\theta_c=\frac{\pi}{2(\log c+1)}$.
We divide the integral in \eqref{eq3.1} into the major arc $|y|\le\theta_cx$ and the minor arcs $\theta_cx<|y|\le\pi$, writing
\begin{align}\label{eq3.2}
a_c(n)=\mathcal{I}_{1,c}(n)+\mathcal{I}_{2,c}(n),
\end{align}
where
\begin{align}\label{I_1,2}
\mathcal{I}_{1,c}(n)=\frac{e^{nx}}{2\pi}
\int_{-\theta_cx}^{\theta_cx}
F_c(e^{-x+iy})e^{-iny}\diff y\quad\text{and}\quad \mathcal{I}_{2,c}(n)=\frac{e^{nx}}{2\pi}
\int_{\theta_cx<|y|\le\pi}
F_c(e^{-x+iy})e^{-iny}\diff y.
\end{align}
\par In Section \ref{subsec3.1}, we derive a positive lower bound for $\mathcal{I}_{1,c}(n)$ using Euler–Maclaurin summation and a saddle-point analysis. Section \ref{subsec3.2} provides an upper bound for $\vert \mathcal{I}_{2,c}(n)\vert $ through uniform estimates near roots of unity. In Section \ref{subsec3.3}, we compare these bounds to establish the positivity of the coefficients $a_c(n)$ for all sufficiently large $n$.
\subsection{Major Arc Analysis}\label{subsec3.1}
Recall from \eqref{gf-Fc} that 
\begin{align*}
F_c(q)=\frac{(-cq;q)_{\infty}}{(-q;q)_{\infty}}=\prod_{m=1}^{\infty}\frac{1+cq^m}{1+q^m}.
\end{align*}
Define
\begin{align*}
f_c(z):=\operatorname{Log}(1+ce^{-z})-\operatorname{Log}(1+e^{-z}),   
\end{align*}
where $\operatorname{Log}(z)$ denotes the principal logarithm for complex $z$. 
The logarithms $\operatorname{Log}(1+ce^{-z})$ and $\operatorname{Log}(1+e^{-z})$
are analytic in the sector
$\mathrm{Re}(z)\ge0$, $|\mathrm{Im}(z)|\le\theta_c\mathrm{Re}(z)$.
Indeed, for $a=1,c$, the real part of $1+ae^{-z}$ is positive:
if $0\le\mathrm{Re}(z)\le\log c+1$, then $|\mathrm{Im}(z)|\le\pi/2$; if
$\mathrm{Re}(z)\ge\log c+1$, then $a|e^{-z}|\le e^{-1}<1$.
Thus, putting $q=e^{-z}$, the normally convergent sum
\begin{align}\label{eq2}
\log F_c(e^{-z})=\sum_{m=1}^{\infty}f_c(mz)
\end{align}
defines an analytic logarithm of $F_c(e^{-z})$, chosen to be real for positive $z$.
The function $f_c$ and its first two derivatives decay exponentially as $\mathrm{Re}(z)\to\infty$ in this sector. 
Applying \eqref{eq-g-EM-sum} to $g(u)=f_c(zu)$, we have
\begin{align}\label{eq3}
\sum_{m=1}^{\infty}f_c(mz)=\frac{1}{z}\int_{0}^{\infty}f_c(t)\diff t-\frac{f_c(0)}{2}-\frac{1}{12}f_c'(0)z+\mathcal{R}_2(z),
\end{align}
where
\begin{align}\label{eq3a}
|\mathcal{R}_2(z)|\le\frac{|z|^2}{12}\int_0^\infty|f_c''(zu)|\diff u.
\end{align}
We first evaluate the integral in \eqref{eq3}. 
The substitution $u=e^{-t}$ gives
\begin{align}\label{eq4}
A_c:=\int_{0}^{\infty}f_c(t)\diff t=\int_{0}^{1}\frac{\log(1+cu)-\log(1+u)}{u}\diff u=-\text{Li}_2(-c)+\text{Li}_2(-1).
\end{align}
Substituting \eqref{eq4} into \eqref{eq3} along with $f_c(0)=\log((c+1)/2)$ and $f_c'(0)=-\frac{c-1}{2(c+1)}$, we obtain
\begin{align}\label{eq5}
\log F_c(e^{-z})=\frac{A_c}{z}-\frac{1}{2}\log\frac{c+1}{2}+\frac{c-1}{24(c+1)}z+\mathcal{R}_2(z).
\end{align}

By means of \eqref{rel-Li2} and $\text{Li}_2(-1)=-\pi^2/12$,
we have
\begin{align*}
A_c=\frac{\pi^2}{12}+\frac{\log^2(c)}{2}+\text{Li}_2\left(-\frac{1}{c}\right).
\end{align*}
Since $\frac{\diff A_c}{\diff c}=\frac{\log(c+1)}{c}>0$,
we have $A_c<A_{c+1}$.
The alternating series gives
$$-\frac12+\frac1{16}-\frac1{72}<\operatorname{Li}_2(-1/2)<-\frac{1}{2}+\frac1{16}.$$
This, along with the elementary bounds $2/3<\log2<7/10$ and $49/5<\pi^2<10$, yields $0.5875<A_2<0.640834$.
Consequently, $A_c>0$ for $c\ge2$.

\par Next, we bound the remainder. Differentiating $f_c(z)$ twice gives
\begin{align*}
f_c''(z)=\frac{ce^{-z}}{(1+ce^{-z})^2}-\frac{e^{-z}}{(1+e^{-z})^2}.
\end{align*}
Let $t=\mathrm{Re}(zu)$.
We consider the following two cases.

(1)~For $0\le t\le \log c+1$,
we have $|\mathrm{Im}(zu)|\le\theta_c t\le\pi/2$.
Hence, for $a=1,c$, we have
\begin{align*}
|1+ae^{-zu}|^2
=1+2ae^{-t}\cos(\mathrm{Im}(zu))+a^2e^{-2t}
\ge1+a^2e^{-2t}.
\end{align*}
Consequently,
\begin{align}\label{eq-doulbe-deriI}
|f_c''(zu)|\le\frac{ce^{-t}}{1+c^2e^{-2t}}+\frac{e^{-t}}{1+e^{-2t}}.
\end{align}

(2)~For $t\ge\log c+1$,
we know $|1+ae^{-zu}|\geq 1-|ae^{-zu}|=1-ae^{-t}$, where $a=1,c$.
Thus,
\begin{align}\label{eq-doulbe-deri}
|f_c''(zu)|\le\frac{ce^{-t}}{(1-ce^{-t})^2}+\frac{e^{-t}}{(1-e^{-t})^2}.
\end{align}
Let $z=x-iy$. Then $t=xu$ and $\diff u=\diff t/x$.
Integrating the two bounds, \eqref{eq-doulbe-deriI} and \eqref{eq-doulbe-deri},
we obtain
\begin{align}\label{W_c}
 &\quad \int_0^{\infty}|f_c''(zu)|\diff u\\[5pt]
&\leq \frac{1}{x}\left(\int_0^{\log c+1}
\left(\frac{ce^{-t}}{1+c^2e^{-2t}}+\frac{e^{-t}}{1+e^{-2t}}\right)\diff t\right. \left.+\int_{\log c+1}^{\infty} \left(\frac{ce^{-t}}{(1-ce^{-t})^2}+\frac{e^{-t}}{(1-e^{-t})^2}\right)\diff t\right)=\frac{W_c}{x},\notag
\end{align}
where
\begin{align*}
W_c=\arctan(c)+\frac{\pi}{4}-\arctan(e^{-1})-\arctan\left(\frac{1}{ec}\right)+\frac1{e-1}+\frac1{ec-1}.
\end{align*}
To bound $W_c$ independently of $c$, let $h(a)=a/(1-a)-\arctan(a)$ for $0<a<1$. Since $h'(a)>0$,
we have $h(1/(ec))\le h(1/(2e))$ for $c\ge2$. Together with
$\arctan(c)=\pi/2-\arctan(1/c)<\pi/2$, this yields
\begin{align*}
W_c<\frac{3\pi}{4}-\arctan(e^{-1})+\frac1{e-1}+h\left(\frac{1}{2e}\right)<2.630.
\end{align*}
Since $|z|^2=x^2+y^2\leq(1+\theta_c^2)x^2$, it follows from \eqref{eq3a} and \eqref{W_c} that
\begin{align*}
|\mathcal{R}_2(z)|\le \frac{2.630(1+\theta_c^2)}{12}x.
\end{align*}
Define $\varepsilon(z):=\frac{c-1}{24(c+1)}z+\mathcal{R}_2(z)$. Then, using $\theta_c<1$, we conclude that
\begin{align*}
|\varepsilon(z)|\le\left(\frac{c-1}{24(c+1)}\sqrt{1+\theta_c^2}+\frac{2.630(1+\theta_c^2)}{12}\right)x<\frac{x}{2}<x.
\end{align*}
Exponentiating \eqref{eq5} therefore gives the uniform
major arc approximation
\begin{align}\label{eq5.5_}
F_c(e^{-z})=\sqrt{\frac{2}{c+1}}\exp\left(\frac{A_c}{z}+\varepsilon(z)\right),
\qquad |\varepsilon(z)|<x.
\end{align}
Thus, $F_c(e^{-z})$ is approximated explicitly by $\sqrt{\frac{2}{c+1}}e^{A_c/z}$ with a uniformly small multiplicative error $e^{\varepsilon(z)}$. The integral $\mathcal{I}_{1,c}(n)$ given by \eqref{I_1,2} contains $e^{nx}F_c(e^{-x+iy})$. The exponential size is $\exp\left(nx+\frac{A_c}{x}\right)$, and the saddle point of $nx+\frac{A_c}{x}$ occurs at $x=\sqrt{\frac{A_c}{n}}$.
\par Substituting \eqref{eq5.5_} into $\mathcal{I}_{1,c}(n)$ given by \eqref{I_1,2} yields
\begin{align}\label{eq3.4}
\mathcal{I}_{1,c}(n)&=\frac{1}{\pi\sqrt{2(c+1)}}\int_{-\theta_cx}^{\theta_cx}\exp\left(\frac{A_c}{z}+nz\right)e^{\varepsilon(z)} \diff y.
\end{align}
Set $y=x^{3/2}t/\sqrt{A_c}$. (This scale makes the quadratic decay in the exponent independent of $x$ to leading order.) Then $\diff y=x^{3/2}\diff t/\sqrt{A_c}$, $|t|\le\theta_c\sqrt{A_c/x}=:J$,
\begin{align*}
z=x\left(1-i\sqrt{\frac{x}{A_c}}t\right),~~
\frac{A_c}{z}=\frac{A_c}{x}\frac{1+i\sqrt{x/A_c}t}{1+(x/A_c)t^2},~~ \frac{A_c}{z}+nz=\frac{2A_c}{x}-\frac{t^2}{1+(x/A_c)t^2}-i\frac{\sqrt{x/A_c}t^3}{1+(x/A_c)t^2},
\end{align*}
at $x=\sqrt{\frac{A_c}{n}}$.
With this substitution, \eqref{eq3.4} becomes
\begin{align}\label{eq3.5}
\mathcal{I}_{1,c}(n)=\frac{1}{\pi\sqrt{2(c+1)A_c}}\,x^{3/2}e^{2A_c/x}
\int_{-J}^{J}G_x(t)e^{\varepsilon(z)}\diff t,
\end{align}
where
\begin{align*}
G_x(t)=\exp\left(-\frac{t^2}{1+(x/A_c)t^2}\right)\exp\left(-i\frac{\sqrt{x/A_c}t^3}{1+(x/A_c)t^2}\right).
\end{align*}
By conjugation symmetry, $\mathcal{I}_{1,c}(n)$ is real. Separating the principal term from the error in \eqref{eq3.5}, we have
\begin{align}\label{eq3.6}
\mathcal{I}_{1,c}(n)=\frac{1}{\pi\sqrt{2(c+1)A_c}}\,x^{3/2}e^{2A_c/x}
\left(\int_{-J}^{J}G_x(t)\diff t+\int_{-J}^{J}G_x(t)\left(e^{\varepsilon(z)}-1\right)\diff t\right).
\end{align}
Assume
\begin{align}\label{eq_x_bd}
x\leq\min\left\{\frac{1}{120},\frac{A_c}{64},\frac{\theta_c^2A_c}{4}\right\}.
\end{align}
Then $J\geq2$. For $|t|\le J$, we have
\begin{align}\label{eq_bdd_g_x}
1+(x/A_c)t^2\leq1+(x/A_c)J^2=1+\theta_c^2\leq2, \qquad \vert G_x(t)\vert\leq e^{-t^2/2}.
\end{align}
Since $\vert \varepsilon(z)\vert\leq x$, the error term satisfies
\begin{align}\label{eq3.7}
\Bigg\vert\int_{-J}^{J}G_x(t)\left(e^{\varepsilon(z)}-1\right)\diff t\Bigg\vert\leq (e^x-1)\int_{\mathbb{R}}e^{-t^2/2}\diff t=\sqrt{2\pi}(e^x-1). 
\end{align}
It remains to bound the first integral in \eqref{eq3.6}
from below. For that, we start with 
\begin{align*}
\mathrm{Re}\left(G_x(t)\right)=\exp\left(-\frac{t^2}{1+(x/A_c)t^2}\right)\cos\left(\frac{\sqrt{x/A_c}t^3}{1+(x/A_c)t^2}\right).
\end{align*}
We consider the following two cases:
(1)~For $|t|\leq2$, condition \eqref{eq_x_bd} gives
\begin{align*}
\Bigg\vert\frac{\sqrt{x/A_c}t^3}{1+(x/A_c)t^2}\Bigg\vert\leq8\sqrt{\frac{x}{A_c}}\leq 1<\frac{\pi}{2},
\end{align*}
and therefore 
\begin{align*}
\cos\left(\frac{\sqrt{x/A_c}t^3}{1+(x/A_c)t^2}\right)>0,
\end{align*}
which further implies that $\mathrm{Re}\left(G_x(t)\right)>0$ for $|t|\leq 2$.
In particular,
for $|t|\leq1$,
\begin{align*}
\Bigg\vert\frac{\sqrt{x/A_c}t^3}{1+(x/A_c)t^2}\Bigg\vert\leq\sqrt{\frac{x}{A_c}}\leq\frac{1}{8}<\frac{\pi}{2}\quad\text{implies}\quad
\cos\left(\frac{\sqrt{x/A_c}t^3}{1+(x/A_c)t^2}\right)\geq \cos\left(1/8\right),
\end{align*}
and
\begin{align*}
\frac{t^2}{1+(x/A_c)t^2}\leq 1\qquad\text{implies}\qquad\exp\left(-\frac{t^2}{1+(x/A_c)t^2}\right)\geq e^{-1}.
\end{align*}
Therefore,
\begin{align}\label{eq_int_G1}
\int_{-1}^{1}\mathrm{Re}\left(G_x(t)\right)\diff t\geq 2e^{-1}\cos\left(1/8\right).
\end{align}

(2)~For $2\leq |t|\leq J$, as in \eqref{eq_bdd_g_x}, we have $|G_x(t)|\leq e^{-t^2/2}$.
The remaining parts, $2\le|t|\le J$, can therefore reduce the integral by at most
\begin{align}\label{eq_int_G2}
2\int_2^{\infty}e^{-t^2/2}\diff t\leq \int_2^{\infty}te^{-t^2/2}\diff t=e^{-2}.
\end{align}
The portions $1<|t|\le2$ are nonnegative and omitting them preserves a lower bound. Consequently, from \eqref{eq_int_G1} and \eqref{eq_int_G2},
\begin{align}\label{eq3.8}
\int_{-J}^{J}\mathrm{Re}\left(G_x(t)\right)\diff t\geq 2e^{-1}\cos\left(1/8\right)-e^{-2}.
\end{align}
Combining \eqref{eq3.6}, \eqref{eq3.7}, and \eqref{eq3.8}, we arrive at the following
\begin{align}\label{eq3.9}
\mathcal{I}_{1,c}(n)\geq\frac{1}{\pi\sqrt{2(c+1)A_c}}\,x^{3/2}e^{2A_c/x}\left(2e^{-1}\cos\left(1/8\right)-e^{-2}-\sqrt{2\pi}(e^x-1)\right).
\end{align}
For $0<x\leq1/120$, the expression in parentheses in \eqref{eq3.9} exceeds $0.5737>1/2$. Thus,  
\begin{align}\label{eq3.10}
\mathcal{I}_{1,c}(n)>\frac{1}{2\pi\sqrt{2(c+1)A_c}}x^{3/2}e^{2A_c/x}.
\end{align}
This establishes the purpose of this section, which was to find a positive lower bound for $\mathcal{I}_{1,c}(n)$.
\subsection{The Minor Arc}\label{subsec3.2}
To estimate the minor arc contribution $\mathcal{I}_{2,c}(n)$, we decompose $F_c(q)$ into three factors whose sizes can be controlled separately.
Substituting $z=-cq$ into \eqref{JTPI}, we get
\begin{align*}
\Theta_c(q):=\sum_{m=-\infty}^\infty c^m q^{m(m+1)/2}
=(-cq;q)_\infty(-1/c;q)_\infty(q;q)_\infty. 
\end{align*}
Since $(q;q)_\infty(-q;q)_\infty=(q^2;q^2)_\infty$, it follows that
\begin{align*}
F_c(q)=\frac{(-cq;q)_\infty}{(-q;q)_\infty}=\frac{\Theta_c(q)}{(q^2;q^2)_\infty(-1/c;q)_\infty}=\Theta_c(q)(q^2;q^2)_{\infty}^{-1}e^{S_c(q)},
\end{align*}
where
\begin{align}\label{eq5.12a}
S_c(q)=-\sum_{m=0}^{\infty}
\operatorname{Log}\left(1+\frac{q^m}{c}\right)\quad\text{and}\quad\mathrm{Re}(S_c(q))=-\log\vert(-1/c;q)_{\infty}\vert.
\end{align}
Here, the logarithm of the infinite product is the analytic branch specified by the convergent sum. Taking absolute values and logarithms gives
\begin{align}\label{eq5.18*}
\log|F_c(q)|= \log|\Theta_c(q)|+\log|(q^2;q^2)_{\infty}|^{-1}+\mathrm{Re}(S_c(q)).
\end{align}
\par In Sections \ref{sec3.2.1}, \ref{sec3.2.2}, and \ref{sec3.2.3}, we find estimates for $\log|\Theta_c(q)|$, $\log|(q^2;q^2)_{\infty}|^{-1}$, and $\mathrm{Re}(S_c(q))$, respectively. These estimates are combined in Section \ref{sec3.2.4} to obtain a uniform upper bound for $|F_c(q)|$, which is then used in Section \ref{sec3.2.5} to bound $|\mathcal{I}_{2,c}(n)|$.
\par Throughout the minor-arc argument, assume $0<x\le 1/120$ and $\theta_c x<|y|\le\pi$. We express $q=e^{-x+iy}$ relative to a nearby root of unity by writing
\begin{align*}
q=\mu e^{-w},\qquad
\mu=e^{2\pi ih/k},\qquad w=x-i\beta.
\end{align*}
For any $y$,
the Dirichlet approximation supplies coprime $h\in\mathbb Z$ and $k\ge1$ with
\begin{align*}
1\le k\le x^{-1/2},\qquad y=\frac{2\pi h}{k}+\beta, \qquad |\beta|\le\frac{2\pi\sqrt{x}}k.
\end{align*}
(Let $N=\lfloor x^{-1/2}\rfloor$ and $\alpha=y/(2\pi)$. By Dirichlet approximation, there exist coprime $h,k$ with $1\leq k\leq N$ and $|k\alpha-h|<1/N$. Among the $N+1$ fractional parts of $0,\alpha,\ldots,N\alpha$, a circular gap is at most $1/(N+1)$. Subtracting the two corresponding multiples gives $1\le k\le N$ and an integer $h$ with $|k\alpha-h|\le1/(N+1)<\sqrt{x}$. This implies that $|y-2\pi h/k|\leq 2\pi\sqrt{x}/k$.)
Set
\begin{align*}
v=\frac\beta x, \quad
\epsilon=\begin{cases}0,&k\text{ odd},\\1,&k\text{ even}.\end{cases}
\end{align*}
Then $w=x(1-iv), |w|^2=x^2+\beta^2=x^2(1+v^2)$ and
\begin{align}
k^2|w|^2=k^2x^2+k^2\beta^2\leq (1+4\pi^2)x. \label{eq-bound-kw}
\end{align}
These parameters will be used throughout the estimates below.
\subsubsection{\underline{Estimate for $\log|\Theta_c(q)|$}}\label{sec3.2.1}
We have
\begin{align}\label{eq4.1}
\Theta_c(q)=\sum_{m=-\infty}^\infty c^m q^{m(m+1)/2}. 
\end{align}
Put $b_m:=\mu^{m(m+1)/2}$. 
Using $m(m+1)=(m+\tfrac12)^2-\tfrac14$, we rewrite
\eqref{eq4.1} as
\begin{align}\label{eq4.2}
\Theta_c(\mu e^{-w})=e^{w/8-\log c/2} \sum_{m=-\infty}^\infty b_m \exp\left(-\frac{w}{2}\left(m+\frac{1}{2}\right)^2+\log c\left(m+\frac{1}{2}\right)\right).
\end{align}
\par We first express the root-of-unity factor $b_m$ as a finite Fourier sum.
\begin{lemma}\label{lemma3.1}
For $b_m=\mu^{m(m+1)/2}=e^{\pi ihm(m+1)/k}$, we have
\begin{align}\label{b_m}
b_m=\sum_{r=0}^{k-1}\gamma_r e^{2\pi i\nu_r m},
\end{align}
where $\nu_r=\frac{r+\epsilon/2}{k}$ and $\gamma_r$ satisfies $|\gamma_r|=\frac{1}{\sqrt{k}}$.
\end{lemma}
\begin{proof}
Since $b_{m+k}=(-1)^{h(k+1)}b_m=(-1)^\epsilon b_m$, the sequence $d_m:=b_m e^{-\pi i\epsilon m/k}$ is $k$-periodic.
Let
\begin{align}\label{gamma_r}
\gamma_r:=\frac1k\sum_{j=0}^{k-1}d_j e^{-2\pi irj/k},\qquad 0\le r<k.
\end{align}
By Fourier inversion directly, we have
\begin{align}
\sum_{r=0}^{k-1}\gamma_r e^{2\pi irm/k}
=\frac1k\sum_{j=0}^{k-1}d_j\sum_{r=0}^{k-1}e^{2\pi ir(m-j)/k}=d_m
=b_m e^{-\pi i\epsilon m/k},\label{bm-eq}
\end{align}
where we use the fact that the inner sum is $k$ if $m\equiv j\pmod{k}$ and zero otherwise, by the finite geometric-series identity. 
Then, \eqref{bm-eq} implies \eqref{b_m}.
\par Now, we prove that $|\gamma_r|=\frac{1}{\sqrt{k}}$. From \eqref{gamma_r}, we have
\begin{align*} 
k^2|\gamma_r|^2
=\sum_{j=0}^{k-1}\sum_{\ell=0}^{k-1}d_j\overline{d_\ell}e^{-2\pi ir(j-\ell)/k}=\sum_{t=0}^{k-1}e^{-2\pi irt/k}\sum_{\ell=0}^{k-1}d_{\ell+t}\overline{d_\ell},
\end{align*}
where we group the pairs by $t\equiv j-\ell\pmod k$ and replace $d_j$ by $d_{\ell+t}$ (as $d_m$ is $k$-periodic). For each $t$, we have
\begin{align*}
\sum_{\ell=0}^{k-1}d_{\ell+t}\overline{d_\ell}=e^{\pi i(ht(t+1)-\epsilon t)/k}\sum_{\ell=0}^{k-1}e^{2\pi iht\ell/k}=\begin{cases}k,&t=0,\\0,&1\leq t\leq k-1.\end{cases}
\end{align*}
Thus $k^2|\gamma_r|^2=k$, proving $|\gamma_r|=\frac{1}{\sqrt{k}}$.
\end{proof}
By using the expression for $b_m$ from Lemma \ref{lemma3.1} in \eqref{eq4.2}, we get
\begin{align}\label{eq4.3}
\Theta_c(\mu e^{-w})=e^{w/8-\log c/2}\sum_{r=0}^{k-1}\gamma_r \mathcal{H}_{r}(w),
\end{align}
where 
\begin{align}\label{eq4.4}
\mathcal{H}_{r}(w)=e^{-\pi i\nu_r}\sum_{m=-\infty}^\infty  \exp\left(-\frac{w}{2}\left(m+\frac{1}{2}\right)^2+\log c\left(m+\frac{1}{2}\right)\right)\exp\left(2\pi i\nu_r\left(m+\frac{1}{2}\right)\right).
\end{align}
\par The next lemma applies Poisson summation to each $\mathcal{H}_{r}(w)$ enabling us to understand the size of the theta factor near each root of unity through explicit Gaussian estimates.
\begin{lemma}\label{lemma3.2}
For $\mathrm{Re}(w)>0$, we have
\begin{align}\label{eq4.4a}
\mathcal{H}_{r}(w)=e^{-\pi i\nu_r}\sqrt{\frac{2\pi}{w}}
\sum_{\ell=-\infty}^\infty (-1)^\ell
\exp\left(\frac{(\log c+2\pi i(\nu_r-\ell))^2}{2w}\right).
\end{align}
\end{lemma}
\begin{proof}
Define
\begin{align*}
\mathcal{G}(t):=\exp\left(-w\frac{t^2}{2}+t\log c \right),\qquad t\in\mathbb R.
\end{align*}
Then, \eqref{eq4.4} becomes
\begin{align*}
\mathcal{H}_{r}(w)=e^{-\pi i\nu_r}\sum_{m=-\infty}^\infty  \mathcal{G}\left(m+\frac{1}{2}\right) \exp\left(2\pi i\nu_r\left(m+\frac{1}{2}\right)\right).
\end{align*}
We use the Fourier transform convention
\begin{align}\label{eq4.6}
\widehat{\mathcal{G}}(\xi)=\int_{\mathbb R}\mathcal{G}(t)e^{-2\pi i\xi t}\diff t.
\end{align}
Since $\mathrm{Re}(w)>0$, the function $\mathcal{G}$ and all its derivatives decay rapidly
on $\mathbb R$, so Poisson summation applies. 
The shifted form of Poisson summation is given by
\begin{align*}
\sum_{m=-\infty}^\infty \mathcal{G}(m+a)e^{2\pi i\nu(m+a)}
=\sum_{\ell=-\infty}^\infty e^{2\pi ia\ell}\widehat{\mathcal{G}}(\ell-\nu),\qquad 
\text{for any } a,\nu\in\mathbb R.
\end{align*}
Taking $a=1/2$ and $\nu=\nu_r$ in the above identity yields
\begin{align}\label{eq4.6a}
\mathcal{H}_{r}(w)=e^{-\pi i\nu_r}\sum_{\ell=-\infty}^\infty
(-1)^\ell \widehat{\mathcal{G}}(\ell-\nu_r).
\end{align}
In view of \eqref{eq4.6} and the Gaussian integral, that is, for any complex $\alpha,\beta$ with $\mathrm{Re}(\alpha)<0$,
\begin{align*}
\int_{\mathbb R}e^{\alpha t^2+\beta t}\diff t
=\sqrt{\frac\pi{-\alpha}}\,e^{-\beta^2/(4\alpha)},
\end{align*}
we have
\begin{align}\label{eq4.7}
\widehat{\mathcal{G}}(\ell-\nu_r)&=\int_{\mathbb R}\exp\left(\frac{-w}{2}t^2+(\log c
+2\pi i(\nu_r-\ell))t\right)\diff t\nonumber\\
&=\sqrt{\frac{2\pi}{w}}
\exp\left(\frac{(\log c+2\pi i(\nu_r-\ell))^2}{2w}\right).
\end{align}
Finally, substituting \eqref{eq4.7} into \eqref{eq4.6a} gives \eqref{eq4.4a}.
\end{proof}
Substituting $\mathcal{H}_{r}(w)$ from Lemma \ref{lemma3.2} into \eqref{eq4.3}, 
we obtain
\begin{align}\label{eq4.8}
\Theta_c(\mu e^{-w})=e^{w/8-\log c/2}\sqrt{\frac{2\pi}{w}}
\sum_{r=0}^{k-1}\gamma_r  e^{-\pi i\nu_r}
\sum_{\ell=-\infty}^\infty (-1)^\ell
\exp\left(\frac{(\log c+2\pi i(\nu_r-\ell))^2}{2w}\right).
\end{align}
Now, we compute the modulus of each transformed term $\exp\left(\frac{(\log c+ia)^2}{2w}\right)$, for any $a\in\mathbb{R}$.
Recalling that $w=x(1-iv)$, we have
\begin{align*}
\mathrm{Re}\left(\frac{(\log c+ia)^2}{2w}\right)=\frac{(\log c)^2-a^2-2 av\log c}{2x(1+v^2)}=\frac{(\log c)^2}{2x}-\frac{(a+v\log c )^2}{2x(1+v^2)}.
\end{align*}
Set
\begin{align*}
\xi:=\frac\epsilon2+\frac{k v\log c}{2\pi},\qquad B:=\frac{2\pi^2x}{k^2|w|^2}.
\end{align*}
Taking $a=2\pi(\nu_r-\ell)$, we have
\begin{align}\label{eq4.9}
&\quad \Bigg\vert\exp\left(\frac{(\log c+2\pi i(\nu_r-\ell))^2}{2w}\right)\Bigg\vert\notag\\[5pt]
&=\exp\left(\mathrm{Re}\left(\frac{(\log c+2\pi i(\nu_r-\ell))^2}{2w}\right)\right)
=e^{\frac{(\log c)^2}{2x}-B(r-k\ell+\xi)^2}.
\end{align}
Applying $|\gamma_r|=1/\sqrt{k}$ and \eqref{eq4.9} to \eqref{eq4.8} yields 
\begin{align}\label{eq4.10}
|\Theta_c(\mu e^{-w})|&\leq e^{x/8}\sqrt{\frac{2\pi}{ck|w|}}
e^{\frac{(\log c)^2}{2x}}\sum_{r=0}^{k-1}\sum_{\ell=-\infty}^\infty
e^{-B(r-k\ell+\xi)^2}\nonumber\\
&=e^{x/8}\sqrt{\frac{2\pi}{ck|w|}}
e^{\frac{(\log c)^2}{2x}}\sum_{j=-\infty}^\infty e^{-B(j+\xi)^2}.
\end{align}
The last equality follows because $r-k\ell$ runs through every integer
exactly once as $0\le r<k$ and $\ell \in\mathbb Z$ (for each $j\in\mathbb{Z}$, choose its unique residue $r$ modulo $k$ then $\ell=(r-j)/k$).
\par It remains to bound the Gaussian sum in \eqref{eq4.10}. Let 
\begin{align}\label{dist_d}
d=\text{dist}(\xi,\mathbb Z)\in[0,1/2].
\end{align}
Integer shifts and reflection leave the Gaussian sum unchanged, so we may replace $\xi$ by $d$ in \eqref{eq4.10}. Pairing the indices $j$ and $-j$ gives
\begin{align*}
e^{Bd^2}\sum_{j=-\infty}^\infty e^{-B(j+\xi)^2}
=1+2\sum_{j=1}^\infty e^{-Bj^2}\cosh(2Bjd).
\end{align*}
Since $0\le d\le1/2$ and $\cosh(x)$ is increasing on $[0,\infty)$,
we have $\cosh(2Bjd)\le\cosh(Bj)$. Thus
\begin{align*}
e^{Bd^2}\sum_{j=-\infty}^\infty e^{-B(j+\xi)^2}
&\leq 1+2\sum_{j=1}^\infty e^{-Bj^2}\cosh(Bj)
= 1+\sum_{j=1}^\infty e^{-Bj^2}\left(e^{Bj}+e^{-Bj}\right)\\
&= 1+\sum_{j=0}^\infty e^{-Bj(j+1)}+\sum_{j=1}^\infty e^{-Bj(j+1)}\\
&=2\sum_{j=0}^\infty e^{-Bj(j+1)}
\leq2\sum_{j=0}^\infty e^{-2Bj}=\frac{2}{1-e^{-2B}},
\end{align*}
where the last inequality uses $j(j+1)\geq2j$ for $j\geq 0$.
Thus,
\begin{align}\label{eq4.11}
\sum_{j=-\infty}^\infty e^{-B(j+\xi)^2}\leq\frac{2e^{-Bd^2}}{1-e^{-2B}}.
\end{align}
From \eqref{eq4.10} and \eqref{eq4.11}, we arrive at
\begin{align*}
|\Theta_c(\mu e^{-w})|\leq e^{x/8}\sqrt{\frac{2\pi}{ck|w|}}
e^{(\log c)^2/(2x)}\frac{2e^{-Bd^2}}{1-e^{-2B}}.
\end{align*}
Taking the logarithm, we get
\begin{align}\label{eq4.12}
\log|\Theta_c(\mu e^{-w})|\leq \frac{x}{8}+\frac{1}{2}\log\left(\frac{2\pi}{ck|w|}\right)+\frac{(\log c)^2}{2x}-Bd^2+\log\left(\frac{2}{1-e^{-2B}}\right).
\end{align}
Finally, by means of \eqref{eq-bound-kw}, we have
\begin{align}\label{B_bd}
B\geq\frac{2\pi^2}{1+4\pi^2}>0.487, \qquad \log\left(\frac2{1-e^{-2B}}\right)<1.17.
\end{align}
Substituting these estimates along with the definition of $B$ in \eqref{eq4.12}, we conclude that
\begin{align}\label{eq4.13}
\log|\Theta_c(\mu e^{-w})|\leq \frac{x}{8}+\frac{1}{2}\log\left(\frac{2\pi}{ck|w|}\right)+\frac{(\log c)^2}{2x}-\frac{2\pi^2d^2}{k^2x(1+v^2)}+1.17.
\end{align}
The negative $d^2$ term is retained because it gives important suppression at even cusps.
\subsubsection{\underline{Estimate for $\log|(q^2;q^2)_{\infty}|^{-1}$}}\label{sec3.2.2}
Recall that
\begin{align*}
\eta(\tau)=e^{\pi i\tau/12}\prod_{m=1}^\infty
\left(1-e^{2\pi im\tau}\right),\quad \text{where } \mathrm{Im}(\tau)>0.
\end{align*}
Set $\tau=2h/k+iw/\pi$.
Then $e^{2\pi i\tau}=q^2$ and $\mathrm{Im}(\tau)=x/\pi$, 
so
\begin{align*}
\eta(\tau)=e^{\pi i\tau/12}(q^2;q^2)_\infty,
\end{align*}
which gives
\begin{align}\label{eq5.1}
\log|(q^2;q^2)_\infty|^{-1}=-\frac{x}{12}-\log|\eta(\tau)|.
\end{align}
To apply the transformation law at a reduced cusp, set $h'=2h/g$ and $k'=k/g$,
where $g=\gcd(k,2)$.
Then $\tau=h'/k'+iw/\pi$ and $\gcd(h',k')=1$, 
which implies that there exist integers $a,b$ such that $-ah'-bk'=1$. 
Define
\begin{align*}
\gamma:=
\begin{pmatrix}
a&b\\
k'&-h'
\end{pmatrix}
\in\operatorname{SL}_2(\mathbb Z),
\qquad
\tau':=\gamma\tau=\frac{a\tau+b}{k'\tau-h'}. 
\end{align*}
It follows from the $\eta$-transformation law \eqref{eq_eta_TF} that
\begin{align*}
\vert\eta(\tau')\vert=\vert k'\tau-h'\vert^{1/2}\vert \eta(\tau)\vert.
\end{align*}
Since $k'\tau-h'=ik'w/\pi$ and $ah'+bk'=-1$, we have $\tau'=\frac{a}{k'}+\frac{i\pi}{k'^2w}$, 
\begin{align*}
|k'\tau-h'|=\frac{k|w|}{\pi g}
\text{ and } \mathrm{Im}(\tau')
=\frac{\pi x}{k'^2|w|^2}
=\frac{\pi g^2}{k^2x(1+v^2)}.
\end{align*}
So
\begin{align}\label{eq5.2}
-\log \vert \eta(\tau)\vert&=\frac{1}{2}\log\vert k'\tau-h'\vert-\log \vert \eta(\tau')\vert=\frac{1}{2}\log\left(\frac{k|w|}{\pi g}\right)-\log \vert \eta(\tau')\vert
\end{align}
Let $Q=e^{2\pi i\tau'}$. 
By the definition of the $\eta$-function, we obtain
\begin{align}\label{eq5.3}
-\log\vert \eta(\tau')\vert&=\frac{\pi}{12}\mathrm{Im}(\tau')-\log\vert(Q;Q)_{\infty}\vert\nonumber\\
&=\frac{\pi^2 g^2}{12k^2x(1+v^2)}-\log\vert(Q;Q)_{\infty}\vert.
\end{align}
Since
\begin{align*}
|Q|=e^{-2\pi\mathrm{Im}(\tau')}=\exp\left(-\frac{2\pi^2g^2}{k^2x(1+v^2)}\right)
=e^{-g^2B}<1,
\end{align*}
we have
\begin{align*}
E_{\eta}:=\vert\log\vert(Q;Q)_{\infty}\vert\vert
&=\Bigg\vert\sum_{m=1}^\infty \log\vert 1-Q^m\vert \Bigg\vert
\leq \sum_{m=1}^\infty -\log (1-\vert Q\vert^m )\\
&\leq \sum_{m=1}^\infty \sum_{j=1}^\infty \frac{|Q|^{mj}}{j}
\leq \sum_{m=1}^\infty \sum_{j=1}^\infty |Q|^{mj}\\
&\le\sum_{m=1}^\infty \frac{|Q|^m}{1-|Q|^m}\le\frac1{1-|Q|}\sum_{m=1}^\infty |Q|^m=\frac{|Q|}{(1-|Q|)^2}.
\end{align*}
Recall from \eqref{B_bd} that $B>0.487$. For odd $k$,
we have $g=1$, so 
$|Q|=e^{-B}<e^{-0.487}<0.615$ and
\begin{align}\label{eq5.4}
\vert\log\vert(Q;Q)_{\infty}\vert\vert<\frac{0.615}{(0.385)^2}<4.15.
\end{align}
For even $k$, we have $g=2$, so $|Q|=e^{-4B}<e^{-1.948}<0.143$ and
\begin{align}\label{eq5.5}
\vert\log\vert(Q;Q)_{\infty}\vert\vert<\frac{0.143}{(0.857)^2}<0.195.
\end{align}
Combining \eqref{eq5.1}, \eqref{eq5.2}, \eqref{eq5.3}, \eqref{eq5.4}, and \eqref{eq5.5}, we conclude that
\begin{align}\label{eq5.6}
\log|(q^2;q^2)_\infty|^{-1}\le
\frac12\log\frac{k|w|}{\pi g}-\frac{x}{12}
+\frac{\pi^2g^2}{12k^2x(1+v^2)}+E_\eta,
\end{align}
where
\begin{align}\label{eq5.66}
E_{\eta}<\begin{cases}4.15,&g=1,\\0.195,&g=2.\end{cases}
\end{align}
\subsubsection{\underline{Estimate for $\mathrm{Re}(S_c(q))$}}\label{sec3.2.3}
From the definition \eqref{eq5.12a} of $S_c(q)$, we have
\begin{align}\label{eq5.12}
S_c(q)&=-\sum_{m=0}^{\infty}\sum_{j=1}^{\infty}\frac{(-1)^{j+1}}{j}\left(\frac{1}{c}q^m\right)^j=\sum_{j=1}^{\infty}\frac{(-1)^jc^{-j}}{j}\sum_{m=0}^{\infty}q^{mj}\nonumber\\
&=\sum_{j=1}^{\infty}\frac{(-1)^jc^{-j}}{j(1-q^j)}=:S_{c,1}(q)+S_{c,2}(q),
\end{align}
where 
\begin{align*}
S_{c,1}(q):=
\sum_{\substack{j=1\\k\mid j}}^{\infty}
\frac{(-1)^jc^{-j}}{j(1-q^j)},
\qquad\text{and}\qquad
S_{c,2}(q):=
\sum_{\substack{j=1\\k\nmid j}}^{\infty}
\frac{(-1)^jc^{-j}}{j(1-q^j)}.
\end{align*}
\par First, we deal with $S_{c,1}(q)$. Let $j=kr$, for some positive integer $r$. Then, $q^j=(e^{2\pi ih/k}e^{-w})^{kr}=e^{-krw}$, and
\begin{align*}
S_{c,1}(q)=\sum_{r=1}^{\infty}\frac{(-1)^{kr}c^{-kr}}{kr(1-e^{-krw})}.
\end{align*}
For $\mathrm{Re}(z)>0$, let $R(z):=\frac{1}{1-e^{-z}}-\frac{1}{z}$.
Then we have
\begin{align*}
S_{c,1}(q)&=\frac{1}{k^2w}\sum_{r=1}^{\infty}\frac{((-1)^{k}c^{-k})^r}{r^2}+\sum_{r=1}^{\infty}\frac{(-1)^{kr}c^{-kr}}{kr}R(krw)\\
&=\frac{\text{Li}_2((-1)^{k}c^{-k})}{k^2w}+\sum_{r=1}^{\infty}\frac{(-1)^{kr}c^{-kr}}{kr}R(krw).
\end{align*}
Since $w=x(1-iv)$, we have 
$\mathrm{Re}\left(\frac{1}{w}\right)=\frac{1}{x(1+v^2)}$,
\begin{align}\label{eq5.13}
\mathrm{Re}(S_{c,1}(q))\leq\frac{\text{Li}_2((-1)^{k}c^{-k})}{k^2x(1+v^2)}+
\sum_{r=1}^{\infty}\frac{c^{-kr}}{kr}\vert R(krw)\vert.
\end{align}
Next, we derive an upper bound for $|R(krw)|$.
\begin{lemma}\label{lemma3.3} For any $r\geq1$, 
\begin{align}\label{eq5.9}
|R(krw)|\leq 1+\frac{2}{krx}.
\end{align} 
Moreover, $|R(krw)|<1$, for $1\leq r\leq5$.
\end{lemma}
\begin{proof}
Since $|z|\ge \mathrm{Re}(z)$, and $|1-e^{-z}|\ge1-|e^{-z}|=1-e^{-\mathrm{Re}(z)}$, we have
\begin{align}\label{eq5.8}
|R(z)|&\leq\frac1{|1-e^{-z}|}+\frac1{|z|}\leq\frac1{1-e^{-\mathrm{Re}(z)}}+\frac{1}{\mathrm{Re}(z)}\nonumber\\
&=1+\frac{1}{e^{\mathrm{Re}(z)}-1}+\frac{1}{\mathrm{Re}(z)}\leq 1+\frac{2}{\mathrm{Re}(z)},
\end{align}
where the last inequality follows from the fact that $e^t-1\geq t$, for $t>0$. Since $\mathrm{Re}(krw)=krx>0$, from \eqref{eq5.8}, we get \eqref{eq5.9}.
\par Next, we derive a strict upper bound for $\vert R(krw)\vert$ for $1\leq r\leq 5$. 
Recall that $w=x-i\beta$, $k\leq1/\sqrt{x}$, and $|\beta|\leq2\pi\sqrt{x}/k$. Thus, $k^2|w|^2=k^2x^2+k^2\beta^2\leq(1+4\pi^2)x$. For $1\leq r\leq 5$, and $x\leq1/120$,
\begin{align}\label{bound-krw}
|krw|\leq5k|w|\leq5\sqrt{(1+4\pi^2)x}\leq5\sqrt{\frac{1+4\pi^2}{120}}<3.
\end{align}
It is known that
\begin{align*}
\frac1{1-e^{-z}}
=\frac12+\frac12\coth\frac z2\qquad\text{and}\qquad\coth u=\frac1u+2u\sum_{m=1}^\infty \frac1{u^2+\pi^2m^2},
\end{align*}
so
\begin{align*}
R(z)=\frac12+2z\sum_{m=1}^\infty \frac1{z^2+4\pi^2m^2}.
\end{align*}
For $|z|\le3<2\pi$, we have
\begin{align*}
|z^2+4\pi^2m^2|\ge4\pi^2m^2-|z|^2=4\pi^2m^2\left(1-\frac{|z|^2}{4\pi^2m^2}\right)\ge4\pi^2m^2\left(1-\frac{|z|^2}{4\pi^2}\right).
\end{align*}
Hence
\begin{align*}
|R(z)|\leq\frac12+\frac{2|z|}{4\pi^2(1-|z|^2/(4\pi^2))}\sum_{m=1}^{\infty}\frac1{m^2}=\frac12+\frac{|z|}{12(1-|z|^2/(4\pi^2))},
\end{align*}
where we used $\sum_{m=1}^\infty m^{-2}=\pi^2/6$.
Notice that the expression on the right increases in $|z|$ for $0\le|z|<2\pi$. 
We have
\begin{align}\label{eq5.7}
|R(z)|\leq \frac12+\frac3{12(1-9/(4\pi^2))}\approx0.824<1,
\end{align}
which, together with \eqref{bound-krw}, implies that $|R(krw)|<1$ for $1\leq r\leq 5$.
\end{proof}
By Lemma \ref{lemma3.3}, we obtain
\begin{align}\label{eq5.11}
\sum_{r=1}^{\infty}\frac{c^{-kr}}{kr}\vert R(krw)\vert
&\leq\sum_{r=1}^{5}\frac{c^{-kr}}{kr}+\sum_{r=6}^\infty\frac{c^{-kr}}{kr}\left(1+\frac{2}{krx}\right)=\frac{1}{k}\sum_{r=1}^\infty\frac{(c^{-k})^r}{r}+\frac{2}{k^2x}\sum_{r=6}^\infty\frac{c^{-kr}}{r^2}\nonumber\\
&\leq-\frac{1}{k}\log(1-c^{-k})+\frac{1}{18k^2x}\sum_{r=6}^\infty(c^{-k})^r\nonumber\\
&=-\frac{1}{k}\log(1-c^{-k})+\frac{c^{-6k}}{18k^2x(1-c^{-k})}
=: L_k(c)+\frac{T_k(c)}{x},
\end{align}
where
\begin{align}\label{def-LT}
L_k(c)=-\frac{1}{k}\log(1-c^{-k})\qquad\text{and}\qquad T_k(c)=\frac{c^{-6k}}{18k^2(1-c^{-k})}.
\end{align}
For fixed $c$, both $L_k(c)$ and $T_k(c)$ decrease with positive integer $k$ (the numerators decrease and the relevant positive denominators increase). Thus,
\begin{align}
L_1(c)=\log\left(\frac{c}{c-1}\right),\quad L_2(c)=\frac{1}{2}\log\left(\frac{c^2}{c^2-1}\right),\quad L_k(c)\leq L_3(c)\quad(k\ge3),\label{L-bound}\\
T_1(c)=\frac1{18c^5(c-1)},\quad T_2(c)=\frac1{72c^{10}(c^2-1)},\quad
T_k(c)\le T_3(c)\quad(k\ge3).\label{T-bound}
\end{align}

\par Now, we move to the sum $S_{c,2}(q)$. We have
\begin{align}\label{eq5.11a}
\mathrm{Re}(S_{c,2}(q))\leq 
\sum_{\substack{j=1\\k\nmid j}}^{\infty}\frac{c^{-j}}{j\vert 1-q^j\vert}
=: E^2_{S}(c).
\end{align}

The following lemma gives a bound for $\frac{c^{-j}}{j|1-q^j|}$.
\begin{lemma}\label{lemma3.4}
For every integer $j\geq 1$ with $k\nmid j$, we have
\begin{align*}
\frac{c^{-j}}{j|1-q^j|}\leq
\frac{c^{-j}}j+\frac{c^{-j}}{xj^2}.
\end{align*}
If $k=2$, we have
\begin{align*}
\frac{c^{-1}}{|1-q|}<0.277.   
\end{align*}
If $k\ge3$, we have
\begin{align*}
\frac{c^{-1}}{|1-q|}
+\frac{c^{-2}}{2|1-q^2|}
<\frac{0.178}{\sqrt{x}}.
\end{align*}
\end{lemma}
\begin{proof}
For any $j\geq1$, we have $|1-q^j|\ge1-|q^j|=1-e^{-jx}$.
Thus
\begin{align*}
\frac{c^{-j}}{j|1-q^j|}
\le
\frac{c^{-j}}{j(1-e^{-jx})}
\le
\frac{c^{-j}}j+\frac{c^{-j}}{xj^2},
\end{align*}
where we used the fact that $\frac1{1-e^{-t}}\le1+\frac1t$ for $t>0$. 

\par Now, we prove a stronger bound for $j=1,2$ using angular separation. Assume $j\sqrt{x}<1$.
Since
\begin{align*}
|1-q^j|=\big|1-e^{-xj+ijy}\big|,
\end{align*}
where
\begin{align*}
jy=\frac{2\pi jh}{k}+j\beta,
\end{align*}
we first estimate the distance of this angle from $2\pi\mathbb Z$. Because $\gcd(h,k)=1$, the condition $k\nmid j$ implies $k\nmid jh$. Thus the unperturbed angle satisfies
\begin{align*}
\text{dist}\!\left(\frac{2\pi jh}{k},2\pi\mathbb Z\right)
\ge\frac{2\pi}{k}.
\end{align*}
Distance to a fixed set changes by at most the size of the perturbation. Consequently,
\begin{align}\label{eq_dist_jy}
\delta:=
\text{dist}(jy,2\pi\mathbb Z)&\ge\frac{2\pi}{k}-j|\beta|\geq \frac{2\pi}{k}-\frac{2\pi j\sqrt{x}}{k}=\frac{2\pi}{k}(1-j\sqrt{x}).
\end{align}
Then, 
\begin{align*}
|1-e^{-xj}e^{ijy}|&=\sqrt{(1-e^{-xj})^2+4e^{-xj}\sin^2(jy/2)}\\
&\ge2 e^{-xj/2}\,\sin(\delta/2)
\ge\frac{2e^{{-xj}/{2}}}{\pi}\delta
=\frac{4}{k} e^{-jx/2}(1-j\sqrt{x}),
\end{align*}
where the last inequality follows from the fact that $\sin t\ge\frac{2t}{\pi}$ 
for $t\in [0,\pi/2]$.
Therefore, 
\begin{align}\label{eq5.14}
\frac{c^{-j}}{j|1-q^j|}
\le
\frac{k\,c^{-j}e^{jx/2}}{4j(1-j\sqrt{x})}.
\end{align}
Now, we focus on the cases $j=1,2$.
For $k=2$, we have $kc^{-1}\leq1$.
Since $x\le1/120$, we have $\sqrt{x}<0.092$, $e^{x/2}<1.005$, and $e^x<1.009$. Therefore, for $c\geq2$, we have 
\begin{align*}
\frac{c^{-1}}{|1-q|}
\le
\frac{e^{x/2}}{4(1-\sqrt{x})}
<
\frac{1.005}{4(1-0.092)}
<0.277.
\end{align*}
For $k\geq3$, 
\begin{align*}
\frac{c^{-1}}{|1-q|}
+\frac{c^{-2}}{2|1-q^2|}&\leq
k\left(\frac{e^{x/2}}{8(1-\sqrt{x})}+\frac{e^x}{32(1-2\sqrt{x})}\right)\\
&<k\left(\frac{1.005}{8(0.908)}+\frac{1.009}{32(0.816)}\right)
<0.178k\leq\frac{0.178}{\sqrt{x}},
\end{align*}
because $k\leq1/\sqrt{x}.$
\end{proof}
Lemma \ref{lemma3.4} forces us to consider cases based on $k$.\\
\textbf{Case I: $k=1$.} There are no indices satisfying $1\nmid j$, so $E^2_{S}(c)=0$.\\
\textbf{Case II: $k=2$.} In this case, the indices $j$ are precisely the positive odd integers. After treating $j=1$ separately, Lemma \ref{lemma3.4} gives
\begin{align*}
E^2_{S}(c)<0.277+\sum_{\substack{j=3\\j\text{ odd}}}^{\infty}\frac{c^{-j}}j+\frac1x\sum_{\substack{j=3\\j\text{ odd}}}^{\infty}\frac{c^{-j}}{j^2}.
\end{align*}
The odd part of the logarithmic series yields
\begin{align*}
\sum_{\substack{j=3\\j\text{ odd}}}^{\infty}\frac{c^{-j}}j=\frac{1}{2}\log\left(\frac{c+1}{c-1}\right)-\frac{1}{c}.
\end{align*}
So we have
\begin{align*}
D_2(c):=\sum_{\substack{j=3\\j\text{ odd}}}^{\infty}\frac{c^{-j}}{j^2}=\frac{\text{Li}_2(1/c)-\text{Li}_2(-1/c)}{2}-\frac1c.
\end{align*}
Therefore, for $k=2$,
\begin{align}\label{eq5.15}
E^2_{S}(c)<0.277+\frac{1}{2}\log\left(\frac{c+1}{c-1}\right)-\frac{1}{c}+\frac{D_2(c)}{x}.
\end{align}
\textbf{Case III: $k\geq3$.}
After treating $j=1,2$ separately, Lemma \ref{lemma3.4} gives
\begin{align*}
E^2_{S}(c)<\frac{0.178}{\sqrt{x}}+\sum_{\substack{j=3\\k\nmid j}}^{\infty}
\left(\frac{c^{-j}}j+\frac{c^{-j}}{xj^2}\right)<\frac{0.178}{\sqrt{x}}+\sum_{j=3}^{\infty}\frac{c^{-j}}j+\frac1x\sum_{j=3}^{\infty}\frac{c^{-j}}{j^2}.
\end{align*}
Here,
\begin{align*}
\sum_{j=3}^{\infty}\frac{c^{-j}}j=-\log(1-c^{-1})-\frac1c-\frac{1}{2c^2},
\end{align*}
and we put
\begin{align*}
D_3(c):=\sum_{j=3}^{\infty}\frac{c^{-j}}{j^2}=\operatorname{Li}_2(1/c)-\frac1c-\frac1{4c^2}.
\end{align*}
Therefore, for $k\geq3$,
\begin{align}\label{eq5.16}
E^2_{S}(c)<\frac{0.178}{\sqrt{x}}-\log(1-c^{-1})-\frac1c-\frac{1}{2c^2}+\frac{D_3(c)}{x}.
\end{align}
\par Finally, combining \eqref{eq5.12}, \eqref{eq5.13}, and \eqref{eq5.11a}, we obtain
\begin{align}\label{eq5.17}
\mathrm{Re}(S_c(q))\leq \frac{\text{Li}_2((-1)^{k}c^{-k})}{k^2x(1+v^2)}+
L_k(c)+\frac{T_k(c)}{x}+E^2_{S}(c),
\end{align}
where $L_k(c)$ and $T_k(c)$ are given by \eqref{def-LT} and the upper bounds on $E^2_{S}(c)$ are given by \eqref{eq5.15} and \eqref{eq5.16}.
\subsubsection{\underline{A uniform bound for $F_c(q)$}}\label{sec3.2.4}
Recall \eqref{eq5.18*}
\begin{align}\label{eq5.18}
\log|F_c(q)|= \log|\Theta_c(q)|+\log|(q^2;q^2)_{\infty}|^{-1}+\mathrm{Re}(S_c(q)).
\end{align}
Substituting \eqref{eq4.13}, \eqref{eq5.6} and \eqref{eq5.17} into \eqref{eq5.18} gives
\begin{align}\label{eq5.19}
\log|F_c(q)|\leq\frac{\Psi_{k,c}(v)}{x}+\frac{x}{24}+\frac{1}{2}\log\left(\frac{2}{cg}\right)+1.17+E_\eta+L_k(c)+\frac{T_k(c)}{x}+E^2_{S}(c),
\end{align}
where
\begin{align}\label{eq5.20}
\Psi_{k,c}(v)=\frac{(\log c)^2}{2}+\frac{\text{Li}_2((-1)^{k}c^{-k})+\pi^2g^2/12-2\pi^2d^2}{k^2(1+v^2)}.
\end{align}
The following lemma bounds the coefficient $\Psi_{k,c}(v)$ of $1/x$.
\begin{lemma}\label{lemma3.5}
For even $k\geq2$, define $P_c=\pi^2/12+\operatorname{Li}_2\left(-1/c\right)$,
\begin{align*}
\sigma_{k,c}=\frac{(\log c)^2}{2}+\frac{\operatorname{Li}_2(c^{-k})-\pi^2/6}{k^2},
 \qquad \phi_{k,c}=\frac{\pi\log c}{k},\qquad
 \alpha_{k,c}=\frac{\sigma_{k,c}+\sqrt{\sigma_{k,c}^2+\phi_{k,c}^2}}{2},\\[5pt]
\beta_c=\frac{(\log c)^2/2+\sqrt{(\log c)^4/4+(\pi\log c/6)^2}}{2},
 \quad \Omega_c=\max\left\{\frac{(\log c)^2}{2}+\frac{\pi^2}{108},\alpha_{4,c},\beta_c\right\}.    
\end{align*}
Then, on the minor arc,
\begin{align*}
\Psi_{1,c}(v)\leq A_c-\frac{P_c\theta_c^2}{1+\theta_c^2},
\qquad \Psi_{2,c}(v)\leq\alpha_{2,c},
\qquad \Psi_{k,c}(v)\leq\Omega_c\quad(k\geq3).   
\end{align*}
\end{lemma}
\begin{proof}
For odd $k$, we have $g=1$ and $(-1)^kc^{-k}=-c^{-k}$.
Therefore,
\begin{align}\label{eq5.21}
\Psi_{k,c}(v)=\frac{(\log c)^2}{2}+\frac{\text{Li}_2(-c^{-k})+\pi^2/12-2\pi^2d^2}{k^2(1+v^2)}\le
\frac{(\log c)^2}{2}+\frac{\text{Li}_2(-c^{-k})+\pi^2/12}
{k^2(1+v^2)}.
\end{align}
The numerator, $\text{Li}_2(-c^{-k})+\pi^2/12$, is positive (for $0<t=c^{-k}\le1/2$, $-t\le\text{Li}_2(-t)<0$), and therefore we may replace $1+v^2$ by $1$. Thus, for odd $k\geq3$,
\begin{align}\label{eqq1}
\Psi_{k,c}(v)\leq \frac{(\log c)^2}{2}+\frac{\text{Li}_2(-c^{-k})+\pi^2/12}{k^2}\le\frac{(\log c)^2}{2}+\frac{\pi^2}{108}.
\end{align}
Substituting $A_c=P_c+(\log c)^2/2$ and $k=1$ into \eqref{eq5.21}, we get
\begin{align}\label{eq5.22}
\Psi_{1,c}(v)\leq  \frac{(\log c)^2}{2}+\frac{P_c}{1+v^2}=A_c-P_c+\frac{P_c}{1+v^2}=A_c-\frac{P_c v^2}{1+v^2}.
\end{align}
On the minor arc, $\theta_cx<|y|\le\pi$, for $k=1$, we have $h=0$ (otherwise $|y-2\pi h|\geq\pi$, contradicting $|\beta|\leq2\pi\sqrt{x}<\pi$). Thus $y=\beta$, and $|v|=|y|/x>\theta_c$. Therefore $1+v^2>1+\theta^2_c$, and \eqref{eq5.22} proves the asserted bound for $k=1$.
\par For even $k$, $g=2$, $\epsilon=1$, and
\begin{align}\label{eq5.23}
\Psi_{k,c}(v)=\frac{(\log c)^2}{2}+\frac{\text{Li}_2(c^{-k})+\pi^2/3-2\pi^2d^2}
{k^2(1+v^2)}.
\end{align}

\textbf{Case I: $0\le |v|\le \pi/(k\log c)$.} In this case, $\frac{1}{2}\leq\frac12+\frac{k|v|\log c}{2\pi}\leq1$. Therefore, \eqref{dist_d} with $\xi=\frac{\epsilon}{2}+\frac{kv\log c}{2\pi}$ and $\epsilon=1$ yields
\begin{align*}
d=\text{dist}\left(
\frac12+\frac{kv\log c}{2\pi},\mathbb Z
\right)=\frac12-\frac{k|v|\log c}{2\pi}.
\end{align*}
Equivalently, $2\pi d=\pi-k|v|\log c$ or $2\pi^2d^2=(\pi-k|v|\log c)^2/2$.
Substituting this into \eqref{eq5.23}, we obtain
\begin{align}\label{eq5.24}
\Psi_{k,c}(v)=\frac{(\log c)^2}{2}+\frac{\bigl(\text{Li}_2(c^{-k})-\pi^2/6\bigr)/k^2+(\pi\log c/k)|v|-(\log c)^2|v|^2/2}{1+|v|^2}=\frac{\sigma_{k,c}+\phi_{k,c} |v|}{1+|v|^2}.
\end{align}
By the definition of $\alpha_{k,c}>0$, we have
\begin{align}\label{sigma_phi}
\alpha_{k,c}^2-\sigma_{k,c}\alpha_{k,c}=\frac{\phi_{k,c}^2}{4},
\qquad
\alpha_{k,c}-\sigma_{k,c}=\frac{\phi_{k,c}^2}{4\alpha_{k,c}}.
\end{align}
For any $u\geq 0$, we have
\begin{align}\label{ineq-u}
\alpha_{k,c}(1+u^2)-(\sigma_{k,c}+\phi_{k,c}u)=
\alpha_{k,c}\left(u-\frac{\phi_{k,c}}{2\alpha_{k,c}}\right)^2\geq0.
\end{align}
Thus, for every $0\leq |v|\leq \pi/(k\log c)$,
\begin{align}\label{eq5.25}
\frac{\sigma_{k,c}+\phi_{k,c}|v|}{1+|v|^2}\leq \alpha_{k,c}.
\end{align}
From \eqref{eq5.24} and \eqref{eq5.25}, we get that $\Psi_{k,c}(v)\leq \alpha_{k,c}$.\\

\textbf{Case II: $|v|\ge \pi/(k\log c)$.} 
In this case, we have
\begin{align*}
\Psi_{k,c}(v)\le\frac{(\log c)^2}{2}+\frac{\text{Li}_2(c^{-k})+\pi^2/3}{k^2(1+|v|^2)}=:H_{k,c}(|v|).
\end{align*}
The function $H_{k,c}$ is decreasing on $[0,\infty)$. At $|v|=\pi/(k\log c)$ we have $d=0$, so \eqref{eq5.24}--\eqref{eq5.25} give
\begin{align*}
\Psi_{k,c}(v)\le H_{k,c}(|v|)\le H_{k,c}(\pi/(k\log c))
=\frac{\sigma_{k,c}+\phi_{k,c}\pi/(k\log c)}{1+(\pi/(k\log c))^2}\le \alpha_{k,c},
\end{align*}
where the last inequality follows from \eqref{ineq-u}.

Combining the two cases, we obtain that for even values of $k$ and for every real $v$, 
\begin{align}\label{eqq2}
\Psi_{k,c}(v)\le \alpha_{k,c}.
\end{align}

For even $k\geq6$, note that
\begin{align*}
\text{Li}_2(c^{-k})=\sum_{r\geq1}\frac{c^{-kr}}{r^2}< \sum_{r\geq1}\frac{1}{r^2}=\frac{\pi^2}{6},
\end{align*}
and therefore, $\sigma_{k,c}\le(\log c)^2/2$. Also, $\phi_{k,c}\le \pi\log c/6$. 
Since $(\sigma+\sqrt{\sigma^2+\phi^2})/2$ is nondecreasing in each variable for $\phi\geq0$, we obtain that, for all even $k\geq6$,
\begin{align}\label{eqq3}
\alpha_{k,c}\le\frac{(\log c)^2/2+\sqrt{((\log c)^2/2)^2+(\pi\log c/6)^2}}{2}=\beta_c.  
\end{align}
The estimate for $k\geq3$ now follows from \eqref{eqq1}, \eqref{eqq2}, and \eqref{eqq3}.
\end{proof}
For convenience, recall \eqref{eq5.19}:
\begin{align}\label{eq5.26}
\log|F_c(q)|\leq \frac{\Psi_{k,c}(v)}{x}+\frac{x}{24}+\frac{1}{2}\log\left(\frac{2}{cg}\right)+1.17+E_\eta+L_k(c)+\frac{T_k(c)}{x}+E^2_{S}(c),
\end{align}
where the bounds on $\Psi_{k,c}(v)$ are given by Lemma \ref{lemma3.5}, for $E_\eta$ by \eqref{eq5.66}, for $L_k(c), T_k(c)$ by \eqref{def-LT}, for $E^2_{S}(c)$ by \eqref{eq5.15} and \eqref{eq5.16}. Define
\begin{align*}
d_{1,c}&=\frac{P_c \theta_c^2}{1+\theta_c^2}-T_1(c),\\
d_{2,c}&=A_c-\alpha_{2,c}-D_2(c)-T_2(c),\\
d_{3,c}&=A_c-\Omega_{c}-D_3(c)-T_3(c).
\end{align*}
We estimate \eqref{eq5.26} separately in the following three cases.\\
\textbf{Case I: $k=1$.} In this case, $g=1$, 
\begin{align*}
\Psi_{1,c}(v)\le A_c-\frac{P_c \theta_c^2}{1+\theta_c^2},\quad \frac12\log\left(\frac{2}{c}\right)\le0,\quad E_\eta<4.15,\quad L_1(c)\leq\log 2, \quad E^2_{S}(c)=0.
\end{align*}
Substituting these bounds into \eqref{eq5.26} gives
\begin{align*}
\log|F_c(q)|&<\frac{A_c-P_c \theta_c^2/(1+\theta_c^2)+T_1(c)}{x}+\frac{x}{24}+1.17+4.15+\log2\\
&\leq \frac{A_c-d_{1,c}}{x}+6.02,
\end{align*}
since $x\le1/120$. 
The restriction $|y|>\theta_cx$ is essential, as the central arc near $q=1$ produces the main contribution and cannot satisfy the same exponential suppression.\\
\textbf{Case II: $k=2$.}
In this case, $g=2$,
\begin{align*}
&\Psi_{2,c}(v)\le\alpha_{2,c},\quad\frac12\log\left(\frac{1}{c}\right)\le\frac12\log\left(\frac{1}{2}\right),\quad E_\eta<0.195,\\
&L_2(c)\leq\frac{1}{2}\log(4/3),\quad E^2_{\mathrm{S}}(c)\le0.277+\frac12\log\left(\frac{c+1}{c-1}\right)-\frac1c+\frac{D_2(c)}{x},
\end{align*}
where $\frac12\log\left(\frac{c+1}{c-1}\right)-\frac1c\leq\frac12\log3-\frac12$.
Substituting these bounds into \eqref{eq5.26} yields
\begin{align*}
\log|F_c(q)|<&\frac{\alpha_{2,c}+D_2(c)+T_2(c)}{x}+\frac{x}{24}+1.17+0.195+0.277+\frac12\log2-\frac12\\
\leq&\frac{A_c-d_{2,c}}{x}+1.5.
\end{align*}
\textbf{Case III: $k\geq3$.}
In this case, $g=1,2$,
\begin{align*}
&\Psi_{k,c}(v)\le\Omega_{c},\quad\frac12\log\left(\frac{2}{cg}\right)\le0,\quad E_\eta<4.15,\\
L_k(c)+\frac{T_k(c)}{x}\le L_3(c)&+\frac{T_3(c)}{x},\quad E^2_{\mathrm{S}}(c)\le\frac{0.178}{\sqrt{x}}-\log\left(1-c^{-1}\right)-\frac1c-\frac{1}{2c^2}+\frac{D_3(c)}{x},
\end{align*}
where $L_3(c)\leq\frac{1}{3}\log(8/7)$ and $-\log\left(1-c^{-1}\right)-\frac1c-\frac{1}{2c^2}\leq \log2-5/8$.
Substituting these bounds into \eqref{eq5.26} gives
\begin{align*}
\log|F_c(q)|&<\frac{\Omega_c+D_3(c)+T_3(c)}{x}+\frac{0.178}{\sqrt{x}}+\frac{x}{24}+1.17+4.15+\frac{1}{3}\log(8/7)+\log2-5/8\\
&\leq \frac{A_c-d_{3,c}}{x}+\frac{0.178}{\sqrt{x}}+5.44.
\end{align*}
\par Putting all three cases together, we have
\begin{align}\label{F_bd}
|F_c(e^{-x+iy})|\leq\begin{cases}
\exp\left(\frac{A_c-d_{1,c}}{x}+6.02\right), &k=1,~\pi\geq |y|>\theta_c x,\\
\exp\left(\frac{A_c-d_{2,c}}{x}+1.5\right), &k=2,\\
\exp\left(\frac{A_c-d_{3,c}}{x}+\frac{0.178}{\sqrt{x}}+5.44\right), &k\geq3.
\end{cases}
\end{align}
To complete the argument, we show that all three gaps $d_{j,c}$ are positive.
\begin{lemma}\label{lemma_d_j,c}
For every $c\geq2$ and $j\in\{1,2,3\}$, we have $d_{j,c}>0$.
\end{lemma}
\begin{proof}
\textit{Positivity of $d_{1,c}$.} Since $\theta_c<1$ and $P_c\ge \frac{\pi^2}{12}-1/2$,
\begin{align*}
\frac{P_c\theta_c^2}{1+\theta_c^2}
\ge\frac{({\pi^2}/{12}-1/2)\pi^2}{8(\log c+1)^2},\qquad
T_1(c)=\frac{1}{18c^5(c-1)}\le\frac{1}{9c^6}=\frac{e^{-6\log c}}9.
\end{align*}
The function $(\log c+1)^2e^{-6\log c}$ decreases for $\log c\ge\log2$. Hence
\begin{align*}
d_{1,c}(\log c+1)^2\ge\frac{({\pi^2}/{12}-1/2)\pi^2}{8}
-\frac{(1+\log2)^2}{576}>0.
\end{align*}

\medskip
\noindent\textit{Positivity of $d_{2,c}$.} Put $J(c)=(\text{Li}_2(1/c)-\text{Li}_2(-1/c))/2$.
Then, for $c\geq2$, 
\begin{align}\label{J_eq}
J(c)=\sum_{\substack{j\geq1\\ \text{odd}}}\frac{c^{-j}}{j^2}\leq \frac{1}{c}+\frac{1}{9c^3}\cdot\frac{1}{1-1/c^2}=\frac{1}{c}\left(1+\frac{1}{9(c^2-1)}\right)\leq\frac{28}{27c}.
\end{align}
Using \eqref{Li_id} for $z=1/c$, we get
$A_c-\sigma_{2,c}=\pi^2/8-J(c)$. Also, $A_c=(\log c)^2/2+P_c$ and 
$\phi_{2,c}^2/4=\pi^2(\log c)^2/16$. Direct substitution gives
\begin{align}\label{quadraticgap}
G:= A_c^2-\sigma_{2,c} A_c-\frac{\phi_{2,c}^2}{4}
=P_c\left(\frac{\pi^2}{8}-J(c)\right)-\frac{(\log c)^2}{2}J(c).
\end{align}
The positive root of the quadratic $\mathcal{Q}(X)=X^2-\sigma_{2,c}X-\phi_{2,c}^2/4$ is $\alpha_{2,c}$, while the other root is negative. Since $A_c>0$, showing $G>0$ also shows
$A_c>\alpha_{2,c}$. Factorization then gives
\begin{align}\label{Q_fact}
A_c-\alpha_{2,c}=\frac{G}{A_c+\alpha_{2,c}-\sigma_{2,c}}.
\end{align}
\par For $c=2$, using
\begin{align*}
J(2)\leq\frac{14}{27},\quad P_2=\frac{\pi^2}{12}+\text{Li}_2(-1/2)>0.374,\quad \frac{(\log c)^2}{2}J(2)<\frac{7}{27}(\log2)^2, 
\end{align*}
in \eqref{quadraticgap}, we obtain that $G>0.1429$. 
Since $A_2>0$, $G>0$ and $\alpha_{2,2}>\sigma_{2,2}$,
we have $A_2>\alpha_{2,2}$. Furthermore,
\begin{align*}
\sigma_{2,2}=\frac{(\log c)^2}{2}+\frac{\text{Li}_2(2^{-2})-\pi^2/6}{4}
\ge\frac{(\log2)^2}{2}-\frac{\pi^2}{24}>-0.171008.
\end{align*}
So
\begin{align*}
A_2+\alpha_{2,2}-\sigma_{2,2}<2A_2+0.171008<1.46.
\end{align*}
Therefore, \eqref{Q_fact} yields $A_2-\alpha_{2,2}>0.1429/1.46$. Also, 
\begin{align*}
D_2(2)=J(2)-1/2<14/27-1/2<0.018519,\qquad\text{and}\qquad
T_2(2)=\frac{1}{221184}.
\end{align*}
Thus,
\begin{align*}
d_{2,2}=(A_2-\alpha_{2,2})-D_2(2)-T_2(2)>\frac{0.1429}{1.46}-0.018519-\frac{1}{221184}>0.
\end{align*}
\par For $c\geq3$, we have $P_c\geq \frac{\pi^2}{12}-\frac{1}{3}$, $J(c)\leq\frac{28}{27c}$,
and since the function $(\log c)^2e^{-\log c}$ has its global maximum $4/e^2$ at $c=e^2$, 
\begin{align*}
\frac{(\log c)^2}{2}J(c)\leq\frac{14(\log c)^2}{27c}
=\frac{14}{27}(\log c)^2e^{-\log c}\leq \frac{14}{27}\cdot\frac{4}{e^2}=\frac{56}{27e^2}.
\end{align*}
Then, \eqref{quadraticgap} implies that
\begin{align*}
G\geq\left(\frac{\pi^2}{12}-\frac{1}{3}\right)\left(\frac{\pi^2}{8}-\frac{28}{81}\right)-\frac{56}{27e^2}>0.153.
\end{align*}
Since $A_c>0$, $G>0$ and $\alpha_{2,c}>\sigma_{2,c}$,
we have $A_c>\alpha_{2,c}$. Moreover,
\begin{align*}
\sigma_{2,c}=\frac{(\log c)^2}{2}+\frac{\text{Li}_2(c^{-2})-{\pi^2}/{6}}{4}
\ge\frac{(\log c)^2}{2}-\frac{\pi^2}{24}>0,
\end{align*}
hence $A_c+\alpha_{2,c}-\sigma_{2,c}<2A_c$. Therefore, \eqref{Q_fact} yields $A_c-\alpha_{2,c}>0.153/(2A_c)$. On the other hand, \eqref{J_eq} implies
\begin{align*}
D_2(c)<\frac{1}{9c(c^2-1)}\leq\frac{1}{8c^3},\qquad T_2(c)=\frac{1}{72c^{10}(c^2-1)}\leq\frac{1}{64c^{12}}.
\end{align*}
Then,
\begin{align*}
D_2(c)+T_2(c)\leq\frac{1}{c^3}\left(\frac{1}{8}+\frac{1}{3^9\cdot 64}\right)<\frac{0.126}{c^3}.
\end{align*}
Since $2A_c<(\log c)^2+\pi^2/6$ and $((\log c)^2+\pi^2/6)e^{-3\log c}$ decrease for $c\geq 3$,
\begin{align*}
2A_c(D_2(c)+T_2(c))<0.126((\log c)^2+\pi^2/6)e^{-3\log c}
\leq \frac{0.126}{27}((\log3)^2+\pi^2/6)<0.04.
\end{align*}
It follows that
\begin{align*}
d_{2,c}=(A_c-\alpha_{2,c})-(D_2(c)+T_2(c))>\frac{0.153}{2A_c}-\frac{0.04}{2A_c}=\frac{0.113}{2A_c}>0.
\end{align*}

\medskip
\noindent\textit{Positivity of $d_{3,c}$.}
We first establish
\begin{equation}\label{Prough}
\Omega_c\le\frac{(\log c)^2}{2}+\frac{\pi^2/3+\text{Li}_2(c^{-4})}{16}=:\omega_c.
\end{equation}
Clearly, $(\log c)^2/2+\pi^2/108\leq \omega_c$. 
From $\beta_c\geq (\log c)^2/2$ 
and $\beta_c\left(\beta_c-(\log c)^2/2\right)=\pi^2(\log c)^2/144$, we get
\begin{align*}
\beta_c-\frac{(\log c)^2}{2}=\frac{\pi^2(\log c)^2}{144\beta_c}\leq\frac{\pi^2}{72}.
\end{align*}
This further implies that
\begin{align*}
\beta_c\leq \frac{(\log c)^2}{2}+\frac{\pi^2}{72}\leq \omega_c.
\end{align*}
Note that $\alpha_{4,c}$ is the positive root of the quadratic $\mathcal{Q}(X)=X^2-\sigma_{4,c}X-\phi_{4,c}^2/4$ and $\mathcal{Q}((\log c)^2/2)<0$. 
Therefore, $\alpha_{4,c}\ge (\log c)^2/2$. Using this along with $\phi_{4,c}=\pi\log c/4$ in \eqref{sigma_phi} yields $\alpha_{4,c}-\sigma_{4,c}\leq\frac{\pi^2}{32}$. This further implies that
\begin{align*}
\alpha_{4,c}\leq \frac{\pi^2}{32}+\sigma_{4,c}=\omega_c.
\end{align*}
This proves \eqref{Prough}. Finally, the use of $P_c\ge \pi^2/12-1/2$,
\begin{align*}
\text{Li}_2(c^{-4})\leq\sum_{j\geq1}\frac{1}{c^{4j}}&\leq\sum_{j\geq1}\frac{1}{2^{4j}}=\frac{1}{15},\qquad D_3(c)\leq\frac{1}{9}\sum_{j\geq3}\frac{1}{c^{j}}\leq\frac{1}{9}\sum_{j\geq3}\frac{1}{2^{j}}=\frac{1}{36},\\
T_3(c)&=\frac{1}{162c^{15}(c^3-1)}\leq\frac{1}{37158912},
\end{align*}
in the definition of $d_{3,c}$ gives
\begin{align*}
d_{3,c}&\geq (A_c-\omega_{c})-D_3(c)-T_3(c)\geq P_c-\frac{\pi^2/3+\text{Li}_2(c^{-4})}{16}-D_3(c)-T_3(c)\\
&\geq \frac{\pi^2}{12}-\frac{1}{2}-\frac{\pi^2/3+1/15}{16}-\frac{1}{36}-\frac{1}{37158912}>0.084>0.
\end{align*}
\par This completes the proof.
\end{proof}
\subsubsection{\underline{Bounding the Minor Arc Integral $\mathcal{I}_{2,c}(n)$}}\label{sec3.2.5}
Recall from \eqref{I_1,2} that
\begin{align}\label{eq5.32}
\mathcal{I}_{2,c}(n)=\frac{e^{nx}}{2\pi }\int_{\theta_cx<|y|\leq \pi}F_c(e^{-x+iy})e^{-iny} \diff y.
\end{align}
For each point on the minor arc, choose a rational approximation $h/k$ to $y/(2\pi)$ as in Section \ref{subsec3.2}. The corresponding denominator $k$ may vary with $y$. Thus, to apply \eqref{F_bd} to the entire integral, we take the maximum of its
three bounds. More precisely, put
\begin{align*}
B_{1,c}(x)=\exp&\left(-\frac{d_{1,c}}x+6.02\right),\qquad B_{2,c}(x)=\exp\left(-\frac{d_{2,c}}x+1.5\right),\\
 B_{3,c}(x)&=\exp\left(-\frac{d_{3,c}}x+\frac{0.178}{\sqrt{x}}+5.44\right).
\end{align*}
Then, \eqref{F_bd} gives the following estimate
\begin{align*}
|F_c(e^{-x+iy})|\leq e^{A_c/x}\max_{1\leq j\leq3}B_{j,c}(x)
 \qquad(\theta_c x<|y|\leq\pi).   
\end{align*}
Since $|e^{-iny}|=1$ and the total length of the minor arc is at most
$2\pi$, it follows from \eqref{eq5.32} that
\begin{align*}
|\mathcal{I}_{2,c}(n)|\leq\frac{e^{nx}}{2\pi}
 \int_{\theta_c x<|y|\leq\pi}|F_c(e^{-x+iy})|\diff y
 \leq e^{nx+A_c/x}\max_{1\leq j\leq3}B_{j,c}(x).    
\end{align*}
At the saddle point $x=\sqrt{A_c/n}$, we have $nx=A_c/x$.
Consequently,
\begin{align}\label{I_bd}
|\mathcal{I}_{2,c}(n)|\leq e^{2A_c/x}\max_{1\leq j\leq3}B_{j,c}(x).
\end{align}
By Lemma \ref{lemma_d_j,c}, each $d_{j,c}$ is positive. Hence every term in the
maximum contains an exponential saving relative to $e^{2A_c/x}$.
\subsection{An Effective Threshold for Coefficient Positivity}\label{subsec3.3}

Fix an integer $c\geq2$. At the saddle point, set
\begin{align*}
s=\frac1x=\sqrt{\frac n{A_c}},\qquad\Gamma_c=\log\left(2\pi\sqrt{2(c+1)A_c}\right),
\end{align*}
and recall the lower bound for $\mathcal{I}_{1,c}(n)$ obtained in \eqref{eq3.10}
\begin{align}\label{eq5.35}
\mathcal{I}_{1,c}(n)>\frac{1}{2\pi\sqrt{2(c+1)A_c}}x^{3/2}e^{2A_c/x}.
\end{align}
Then, from \eqref{I_bd}, we obtain
\begin{align}\label{I2}
\frac{|\mathcal{I}_{2,c}(n)|}{\mathcal{I}_{1,c}(n)}\le\exp\left(-\min\left\{V_{1,c}(s),V_{2,c}(s),V_{3,c}(s)\right\}\right),
\end{align}
where
\begin{align*}
V_{1,c}(s)&=d_{1,c}s-\frac32\log s-6.02-\Gamma_c,\\
V_{2,c}(s)&=d_{2,c}s-\frac32\log s-1.5-\Gamma_c,\\
V_{3,c}(s)&=d_{3,c}s-0.178\sqrt s-\frac32\log s-5.44-\Gamma_c.
\end{align*}
We choose a threshold beyond which all three functions are increasing. Let
\begin{align*}
\mathfrak{s}_c=\max\left\{120,\frac{64}{A_c},\frac{4}{\theta_c^2A_c},
\frac3{d_{1,c}},\frac3{d_{2,c}},\frac6{d_{3,c}},
\left(\frac{0.356}{d_{3,c}}\right)^2\right\}.
\end{align*}
The first three terms ensure \eqref{eq_x_bd} so that the preceding major-arc and minor-arc estimates apply. The remaining terms ensure that, for $s\ge \mathfrak{s}_c$,
\begin{align*}
V'_{1,c}(s)&\ge d_{1,c}/2>0,\qquad V'_{2,c}(s)\ge d_{2,c}/2>0,\\
V'_{3,c}(s)&=d_{3,c}-0.089/\sqrt s-1.5/s\ge d_{3,c}/2>0.
\end{align*}
Since $d_{j,c}>0$, each $V_{j,c}(s)$ tends to $+\infty$ as
$s\to\infty$. We may therefore define the finite, effectively
computable integer
\begin{equation}\label{Kc}
K_c=\min\left\{j\in\mathbb N:j\ge \mathfrak{s}_c, 
V_{1,c}(j)\ge1,\ V_{2,c}(j)\ge1,\ V_{3,c}(j)\ge1\right\}.
\end{equation}
Set $N_c=\lceil A_cK_c^2\rceil$.
If $n\ge N_c$, then $s=\sqrt{n/A_c}\ge K_c\geq \mathfrak{s}_c$.
By monotonicity and \eqref{Kc}, $V_{j,c}(s)\geq1$ for $j=1,2,3$.
It follows from \eqref{I2} that
\begin{align}\label{eq5.34}
|\mathcal{I}_{2,c}(n)|\leq e^{-1}\mathcal{I}_{1,c}(n).
\end{align}
Finally, from \eqref{eq3.2}, \eqref{eq5.35}, and \eqref{eq5.34}, we arrive at the following
\begin{align*}
a_c(n)&=\mathcal{I}_{1,c}(n)+\mathcal{I}_{2,c}(n)\geq \mathcal{I}_{1,c}(n)-\vert\mathcal{I}_{2,c}(n)\vert\\
&\geq(1-e^{-1})\mathcal{I}_{1,c}(n)>\frac{1-e^{-1}}{2\pi\sqrt{2(c+1)A_c}}x^{3/2}e^{2A_c/x}>0.
\end{align*}
Thus, for all $c\geq2$, $a_c(n)>0$, whenever $n\geq N_c$. This proves Theorem \ref{Thm_main}.
\par For $2\leq c\leq 10$, Table \ref{Table1} lists the values of $A_c$, $K_c$, and $N_c$.
\begin{table}[h]\caption{Values of $A_c$, $K_c$, and $N_c$, for $2\leq c\leq 10$}
\begin{center}
\begin{tabular}{cccc}\label{Table1}
$c$ & $A_c$ (rounded) & $K_c$ & $N_c$ \\
\hline
2 & 0.614279333460 & 120 & 8,846 \\
3 & 1.116908387343 & 120 & 16,084 \\
4 & 1.547472763574 & 120 & 22,284 \\
5 & 1.926812092637 & 120 & 27,747 \\
6 & 2.267475017664 & 120 & 32,652 \\
7 & 2.577695011004 & 120 & 37,119 \\
8 & 2.863208967333 & 125 & 44,738 \\
9 & 3.128196744820 & 132 & 54,506 \\
10 & 3.375810853434 & 138 & 64,289 \\
\end{tabular}
\end{center}
\end{table}
In particular, Theorem \ref{Thm_main} proves positivity for $c=2$ when $n\geq8846$ and for $c=3$ when $n\geq16084$. The remaining finite ranges, $35\leq n\leq8845$ and $28\leq n\leq16083$, respectively, are verified by using \texttt{Mathematica}. This establishes Conjectures \ref{conj1} and \ref{conj2}.

\section{Proof of Theorem \ref{Thm_main_2}}\label{Sec4}
We prove that $q(n)>a_2(n)$ for every integer $n>10$. To do so, we give an upper bound on $a_2(n)$ and a lower bound on $q(n)$ in Propositions \ref{prop1} and \ref{prop2}, respectively.
\begin{lemma}\label{lemma4.1}
For every real $t>0$,
\begin{equation}\label{eqs3}
F_2(e^{-t})<\exp\left(\frac5{8t}\right).
\end{equation}
\end{lemma}
\begin{proof}
Define
\begin{align*}
f_2(s):=\log\frac{1+2e^{-s}}{1+e^{-s}}\qquad(s\ge0). 
\end{align*}
This function is positive and decreasing, since
\begin{align*}
f_2'(s)=-\frac{e^{-s}}{(1+2e^{-s})(1+e^{-s})}<0.  
\end{align*}
Moreover, $f_2(s)=O(e^{-s})$ at infinity,
so $f_2$ is integrable on $[0,\infty)$. On each interval $[(j-1)t,jt]$, monotonicity gives
\begin{align*}
tf_2(jt)\le\int_{(j-1)t}^{jt}f_2(s)\diff s.    
\end{align*}
Summing these inequalities over $j$ and using \eqref{eq2} yields
\begin{align}\label{eqs1}
\log F_2(e^{-t})=\sum_{j=1}^\infty f_2(jt)\le\frac 1t\int_0^\infty f_2(s)\diff s.
\end{align}
With the substitution $u=e^{-s}$, the integral becomes
\begin{align}\label{eqs2}
\int_0^\infty f_2(s)\diff s=\int_0^1\frac{\log(1+2u)-\log(1+u)}u\diff u=A_2<\frac{5}{8}.
\end{align}
Equations \eqref{eqs1} and \eqref{eqs2} prove \eqref{eqs3}.
\end{proof}
\begin{proposition}\label{prop1}
For every integer $n\ge35$,
\begin{equation}\label{eqs5}
 a_2(n)<2e^{\sqrt{5n/2}}.
\end{equation}
\end{proposition}
\begin{proof}
By Corollary \ref{cor1}, $a_2(n)>0$, for all $n\geq35$. The only negative coefficients in $F_2(q)$ up to degree $34$ are $a_2(4)=-1$ and $a_2(13)=a_2(17)=-3$. Therefore, every coefficient of
\begin{align*}
\widehat F_2(q):=F_2(q)+q^4+3q^{13}+3q^{17}
\end{align*}
is nonnegative.
For $n\ge35$, the coefficient of $q^n$ in $\widehat F_2(q)$ is $a_2(n)$.
Evaluating at $q=e^{-t}\in(0,1)$, we may therefore
bound this single term by the full series:
\begin{equation}\label{eqs4}
a_2(n)e^{-nt}\le \widehat F_2(e^{-t})\le F_2(e^{-t})+7,
\end{equation} 
as $e^{-4t}+3e^{-13t}+3e^{-17t}\le7$. By Lemma \ref{lemma4.1} and \eqref{eqs4},
\begin{align*}
a_2(n)<e^{nt}\left(e^{5/(8t)}+7\right)\qquad(t>0). 
\end{align*}
Choose $t=\sqrt{5/(8n)}$. Then
\begin{align*}
a_2(n)<e^{\sqrt{5n/8}}\left(e^{\sqrt{5n/8}}+7\right).
\end{align*}
For $n\ge35$, we have $\sqrt{5n/8}\ge\sqrt{175/8}>4$, and hence $e^{\sqrt{5n/8}}>e^4>7$.
Thus,
\begin{align*}
a_2(n)<2e^{2\sqrt{5n/8}}=2e^{\sqrt{5n/2}},
\end{align*}
as required.
\end{proof}
Now, we shift our attention to $q(n)$. Let $p(n)$ denote the number
of unrestricted partitions of $n$. We use the inequality
\begin{equation}\label{eqs6}
 q(n)\ge p\left(\left\lfloor\frac n2\right\rfloor\right)\qquad(n\geq0).
\end{equation}
This follows from the identities in \cite[Theorem 2.1]{AA_2020}, which were proved
by using Jacobi's identity. For completeness, we give a direct
combinatorial proof.
\par By Euler’s partition theorem, $q(n)$ also counts partitions of $n$ into odd parts.
Write $n=2m+\varepsilon$, where $m=\lfloor n/2\rfloor$ and $\varepsilon\in\{0,1\}$. To a partition
$\lambda=(\lambda_1,\ldots,\lambda_r)$ of $m$, associate
\begin{align*}
\Phi(\lambda)=
 \left(2\lambda_1-1,\ldots,2\lambda_r-1,
       \underbrace{1,\ldots,1}_{r+\varepsilon\text{ additional parts}}\right).
\end{align*}
Every part in $\Phi(\lambda)$ is odd, and the sum of the parts is $\sum_{j=1}^{r}(2\lambda_j-1)+(r+\varepsilon)=2m+\varepsilon=n$. To see that $\Phi$ is injective, let $\ell$ be the number of parts of $\Phi(\lambda)$. Since $\ell=2r+\varepsilon$, the value $r=(\ell-\varepsilon)/2$ is determined by the image. Remove $r+\varepsilon$ parts equal to $1$, and replace each remaining part
$\mu$ by $(\mu+1)/2$. This recovers $\lambda$.
\par We also use the following explicit estimate of Banerjee et al. \cite[Theorem~1.3]{BPRZ}:
\begin{equation}\label{eqs7}
p(m)>\frac{e^{\pi\sqrt{2m/3}}}{4\sqrt3 m}\left(1-\frac1{2\sqrt m}\right)\qquad(m\ge1).
\end{equation}
\begin{proposition}\label{prop2}
For every integer $n\ge2$,
\begin{equation}\label{y1}
 q(n)>\frac1{7n}
 \exp\!\left(\frac\pi{\sqrt3}\sqrt{n-1}\right).
\end{equation}
\end{proposition}
\begin{proof}
Put $m=\lfloor n/2\rfloor\ge1$. Since
$1-1/(2\sqrt m)\ge1/2$, \eqref{eqs6} and
\eqref{eqs7} imply
\begin{align*}
q(n)\ge p(m)>\frac{e^{\pi\sqrt{2m/3}}}{8\sqrt3m}.   
\end{align*}
We use $2m\ge n-1$ in the numerator and $m\le n/2$ in the denominator to obtain
\begin{align*}
 q(n)>\frac{e^{\pi\sqrt{(n-1)/3}}}{4\sqrt3\,n} >\frac{e^{\pi\sqrt{(n-1)/3}}}{7n},
\end{align*}
where the last inequality follows from $4\sqrt3<7$.
\end{proof}
Now, we are ready to prove Theorem \ref{Thm_main_2}.
\begin{proof}[Proof of Theorem \ref{Thm_main_2}]
Using the elementary bounds
\begin{align*}
\frac\pi{\sqrt3}>\frac{181}{100},\qquad
 \sqrt{\frac52}<\frac{791}{500},\qquad
 \sqrt{1-u}\ge1-u\quad(0\le u\le1),  
\end{align*}
for $n\geq 2500$, we have
\begin{align*}
\frac{\pi\sqrt{(n-1)/3}-\sqrt{{5n}/{2}}}{\sqrt n}=\frac\pi{\sqrt3}\sqrt{1-\frac1n}-\sqrt{\frac52}>\frac{181}{100}\left(1-\frac{1}{2500}\right)-\frac{791}{500}
=\frac{56819}{250000}
 >\frac9{40}.
\end{align*}
It follows that for $n\geq 2500$,
\begin{align}\label{eqs9}
 \pi\sqrt{\frac{n-1}{3}}-\sqrt{\frac{5n}{2}}>\frac9{40}\sqrt n.
\end{align}
Thus, from \eqref{eqs5} and \eqref{y1}, we have
\begin{align*}
    q(n)>\frac1{7n}
 \exp\!\left(\pi\sqrt{\frac{n-1}{3}}\right)
 >\frac1{7n}
 \exp\!\left(\sqrt{\frac{5n}{2}}+\frac9{40}\sqrt n\right)
 >2 e^{\sqrt{{5n}/{2}}}>a_2(n),
\end{align*}
where the penultimate inequality uses $9\sqrt{n}/40>\log(14n)$ (for $n\geq 2500$).

Finally, exact coefficient computation by \texttt{Mathematica} verifies $q(n)>a_2(n)$ for
every integer $11\leq n\leq2499$. This completes the proof.
\end{proof}
\section{Concluding Remarks}\label{Sec5}
We have established eventual strict positivity of the coefficients of $F_c(q)=(-cq;q)_\infty/(-q;q)_\infty$ for every integer $c\ge2$, with an effectively computable threshold. The specializations $c=2$ and $c=3$, together with finite verification, settle the conjectured parity inequalities for partitions and overpartitions. We have also proved the conjectured upper bound for the ordinary partition parity difference in terms of the distinct partition function. The explicit estimates in our argument turn the dominance of the contribution near $q=1$ into an effective criterion for positivity, despite the nonmodularity of $F_c(q)$.

A natural direction for further study is to sharpen the major and minor arc estimates, reduce the resulting thresholds, and determine the least positivity threshold for each $c$. It would also be useful to find direct combinatorial proofs of the three inequalities, based on the combinatorial interpretations developed in \cite{AB_2026}. Such proofs could provide a structural explanation for the parity biases and the exceptional small values of $n$.

\bigskip
\noindent \textbf{Acknowledgements.}
This work was supported by
the National Natural Science Foundation of
China (12201177), the Natural Science Foundation of Hebei Province
(A2024205012), Hebei Key Talent Program (Overseas Returnee Platform) (A2026023);
the 111 Center (D26018).

\bigskip

\noindent \textbf{Declaration of AI use.}
ChatGPT (OpenAI) was used to assist in reducing the error term
and reducing the threshold for $n$ in the 
circle-method proof of Theorem \ref{Thm_main} 
in Section \ref{Sec3}.
All AI-assisted calculations, arguments, and suggestions
were independently verified and revised by the authors,
who take full responsibility for the content of the article.


\end{document}